\documentclass{amsart}
\usepackage{amsfonts}
\usepackage{hyperref}
\usepackage{amsmath}
\usepackage{amssymb}
\usepackage{amsthm}
\usepackage{mathrsfs}
\allowdisplaybreaks[4]

\newtheorem{theorem}{Theorem}[section]
\theoremstyle{plain}

\newtheorem{claim}[theorem]{Claim}

\newtheorem{corollary}[theorem]{Corollary}

\newtheorem{definition}[theorem]{Definition}
\newtheorem{example}[theorem]{Example}

\newtheorem{lemma}[theorem]{Lemma}

\newtheorem{proposition}[theorem]{Proposition}
\newtheorem{remark}[theorem]{Remark}

\numberwithin{equation}{section}
\numberwithin{figure}{section}

\begin{document}
\title{\textbf{On a natural $L^2$ metric on the space of Hermitian metrics}}
\author{Jinwei Gao}
\address{Department of Mathematical Sciences, Tsinghua University, Beijing 100084,China. \textit{Email: gaojw19@mails.tsinghua.edu.cn}}

\begin{abstract}
We investigate the space of Hermitian metrics on a fixed complex vector
bundle. This infinite-dimensional space has appeared in the study of
Hermitian-Einstein structures, where a special $L^2$-type Riemannian metric
is introduced. We compute the metric spray, geodesics and curvature associated
to this metric, and show that the exponential map is a diffeomorphism. Though
being geodesically complete, the space of Hermitian metrics is metrically
incomplete, and its metric completion is proved to be the space of
``$L^2$ integrable'' singular Hermitian
metrics. In addition, both the original space and its completion are CAT(0).
In the holomorphic case, it turns out that Griffiths seminegative/semipositive
singular Hermitian metric is always $L^2$ integrable in our sense. Also, in
the Appendix, the Nash-Moser inverse function theorem is utilized to prove that,
for any $L^2$ metric on the space of smooth sections of a given fiber
bundle, the exponential map is always a local diffeomorphism, provided that
each fiber is nonpositively curved.
\end{abstract}
\maketitle
\tableofcontents

\section{Introduction}

Throughout the paper, unless otherwise specified, $M$ is a fixed closed,
oriented, $n$-dimensional smooth manifold, and $E$ is a smooth complex vector
bundle over $M$ of rank r (although in Section 5, we shall further
assume that $M$ is a complex manifold and $E$ is a holomorphic vector bundle).
We also fix a smooth Riemannian metric $g$ on $M$ and denote the volume form
by $\mu_g$. From the Riesz representation theorem (cf. \cite{Ru87},
Chap 2), the positive top form $\mu_g$ induces a Radon measure on $M$, which
we again denote by $\mu_g$. If we say something holds almost everywhere, we
mean that it holds outside of a $\mu_g$-nullset.

Let $Herm(E_x)$ be the real vector space of Hermitian sesquilinear forms on
the fiber $E_x$. The disjoint union of all $Herm(E_x)$ forms a real vector
bundle over $M$, which we denote by $Herm(E)$; in other words, $Herm(E)_x
=Herm(E_x)$. Any smooth local frame $e_1,\dots,e_r$ of $E$ gives rise to
a local trivialization for $Herm(E)$ with the model space $Herm(r)$, by taking
the matrix representation of the Hermitian forms at each point. The space
$\Gamma(Herm(E))$, endowed with smooth compact-open topology, is then a
Fr\'echet space.

By restricting our attention solely to the positive definite ones, we get open
subspaces $Herm^+(E_x)\subseteq Herm(E_x)$ and $Herm^+(E)\subseteq
Herm(E)$. Here $Herm^+(E)$ is actually a fiber bundle over $M$ with fiber
$Herm^+(r)$ and the above local trivialization for $Herm(E)$ restricts to a
local trivialization for $Herm^+(E)$. In this way $Herm^+(E)$ is an open
subbundle of $Herm(E)$. The space of Hermitian metrics on $E$ is just the
space of smooth sections of $Herm^+(E)$, which we denote by $\Gamma
(Herm^+(E))$. It is clear that $\Gamma(Herm^+(E))$ is an open convex cone
in $\Gamma((Herm(E))$.

In his book (\cite{K87}, Chap 6), S. Kobayashi has written down the following natural Riemannian
metric on $\Gamma(Herm^+(E))$ explicitly (Of course it is already implicitly used in \cite{D85} and later in \cite{Siu87}; we shall talk more about it in Remark 2.3):

Given $h\in\Gamma(Herm^+(E))$
and $v,w\in\Gamma(Herm(E))\cong T_h\Gamma(Herm^+(E))$, the inner product
is defined by

\begin{equation}
(v,w)_h\triangleq\int_Mtr\left(h^{-1}vh^{-1}w\right)\mu_g. 
\label{metric}
\end{equation}

Here, $tr(h^{-1}vh^{-1}w)$ is a smooth function on $M$, given in local
coordinates by $h^{\bar{k}i}h^{\bar{j}l}v_{i\bar{j}}w_{l\bar{k}}$. Although $h,v,w$ are all complex objects,
$tr(h^{-1}vh^{-1}w)$ is actually real and $(\cdot,\cdot)_h$ is indeed a real
inner product on $\Gamma(Herm(E))$. We shall explain this metric in detail in
the next section. (Actually, the metric we are dealing with is carried with a
smooth parameter function $\alpha>-1/r$, and \eqref{metric} is the $\alpha=0$ case.)

We remark that this metric is fundamentally different than the Ebin's metric
on the space of all Riemannian metrics on $M$ (cf. \cite{Cla13a}, \cite{FG89},
\cite{GMM91}), since $\mu_g$ here is fixed, while for the Ebin's metric, the
volume element $\mu_g$ is of course varying with the base point $g$.
Actually, the geometry of our metric is much simpler than that of Ebin's metric.

On the other hand, the notion of a singular Hermitian metric on a complex
vector bundle has been introduced in \cite{BP08}, see also \cite{HPS18}, \cite{I20},
\cite{PT18}, \cite{R15}, among others. Roughly speaking, the word ``singular'' here means two things at once:
first, that the metric is not necessarily $C^{\infty}$; second, that certain
vectors in the fibers of the vector bundle are allowed to have either infinite
length or length zero. Usually, one also requires a singular Hermitian
metric to be almost everywhere positive definite, as in \cite{HPS18},
\cite{I20}, \cite{PT18}. Singular Hermitian metrics will be used in our
description of the metric completion of the space of Hermitian metrics.

Our paper is organized as follows. After introducing the background materials
we need in Section 2, we start at Section 3 by computing the metric spray,
geodesics and curvature associated with our $L^2$ metric \eqref{metric}.
Moreover, the corresponding exponential map $Exp_h:\Gamma
(Herm(E))\rightarrow\Gamma(Herm^+(E))$ is shown to be a diffeomorphism.

The first half of Section 4 is devoted to the proof of the following formula
($d$ is the induced distance function on the space of Hermitian metrics):

\begin{equation}
d(h_1,h_2)=\left[\int_Md_x(h_1(x),h_2(x))^2\mu_g(x)\right]^{\frac12}.
\label{equality}
\end{equation}

For the space of Riemannian metrics, Clarke has proved a similar distance
equality (\cite{Cla13b}, Theorem 3.8) for the
Ebin's metric, and characterized the corresponding metric completion in \cite{Cla13a} using a
notion of $\omega$-convergence. His method for proving the distance equality
cannot be applied to our case, for there is no way to obtain a volume-based
estimate on $d$ as in Proposition 4.1 of \cite{Cla13a}, since what we are
dealing with are Hermitian metrics on vector bundles, not Riemannian metrics
on manifolds. This volume-based estimate plays a major role in both his
papers and forms the basis of his proof in \cite{Cla13b}, Theorem 3.8. Instead,
we give a simple proof of \eqref{equality}, using explicit geodesics in $\Gamma(Herm^+(E))$ calculated in Section 3.

With the help of \eqref{equality}, in the rest of Section 4, we derive one of
our main results:

\begin{theorem}
$\Gamma(Herm^+(E))$ is geodesically complete yet metrically incomplete, and
its metric completion can be identified with $L^2(Herm^+(E))$, the space
of $L^2$ integrable singular Hermitian metrics. Both these spaces are CAT(0).
\end{theorem}

Half of our proof is based on the approximation result given in \cite{Ca23},
and we do not use Clarke's more involved notion of ``$\omega
$-convergence''. Again inspired by the recent paper
\cite{Ca23}, an alternative description of the metric completion is given,
although this construction is only valid when the parameter function $\alpha$ is a constant.

In Section 5, we move ahead to examine the $L^2$ integrable condition more
closely and show that, in the holomorphic case, the following holds:

\begin{theorem}
A Griffiths seminegative/semipositive singular Hermitian metric is always
$L^2$ integrable.
\end{theorem}

In addition, the relevant cases of the space of Hermitian metrics on a complex
manifold and the high-order Sobolev metrics are discussed.

Finally, in the Appendix, we explore the exponential map of the general
$L^2$ metric defined on the space of smooth sections of a given fiber
bundle. In particular, using the Nash-Moser inverse function theorem, we show that:

\begin{theorem}
If each fiber is nonpositively curved, then the exponential map for the
general $L^2$ metric is a local diffeomorphism.
\end{theorem}
The last result can be viewed as a version of Cartan-Hadamard theorem in the tame
Fr\'echet setting.

\section{Preliminaries}

In this section, we review some background materials we need.

\subsection{The natural \texorpdfstring{$L^2$}{L2} metric on the space of Hermitian metrics}

Our metric is labeled with a parameter $\alpha$. When $\alpha=0$, it is the
one given in \cite{K87}. For the readers' convenience, we give a detailed
description here.

Let us begin with the space of $r\times r$ Hermitian matrices; it is a real
vector space of dimension $r^2$ and we denote it by $Herm(r)$. The set of
positive-definite $r\times r$ Hermitian matrices, written as $Herm^+(r)$, is
an open convex cone in $Herm(r)$. Therefore $Herm^+(r)$ is trivially an
$r^2$-dimensional smooth manifold whose tangent space at each point is
identified with $Herm(r)$. We fix a real number $\alpha>-1/r$.

There is a natural Riemannian metric on $Herm^+(r)$, defined by:

\begin{equation}
\begin{aligned}
\left\langle v,w\right\rangle_h^{\alpha}  &  \triangleq h^{\bar{k}
i}h^{\bar{j}l}v_{i\bar{j}}w_{l\bar{k}}+\alpha
h^{\bar{k}i}v_{i\bar{k}}h^{\bar{l}j}w_{j\bar{l}}\\
&  =tr(h^{-1}vh^{-1}w)+\alpha tr(h^{-1}v)tr(h^{-1}w),
\label{metric0}
\end{aligned}
\end{equation}
where $h\in Herm^+(r)$ and $v,w\in T_h(Herm^+(r))\cong Herm(r)$.

In our notation, $h^{-1}=(h^{\bar{k}l})$ is the inverse matrix of $h$
so that $\sum h_{i\bar{k}}h^{\bar{k}j}=\delta_i^j=\sum
h^{\bar{j}k}h_{k\bar{i}}$. Note that $h^{-1}$ is also a
positive definite Hermitian matrix.

We must check that $\left\langle\cdot,\cdot\right\rangle_h^{\alpha}$ is
indeed a real inner product on $Herm(r)$ (it has been omitted from
\cite{K87}). First, despite that $h,v,w$ are all complex objects,
$tr(h^{-1}vh^{-1}w)$ and $tr(h^{-1}v)$, $tr(h^{-1}w)$ are actually real:

\begin{align*}
\overline{tr(h^{-1}vh^{-1}w)} & =h^{\bar{i}k}(h^{\bar{l}
j}v_{j\bar{i}})w_{k\bar{l}}=tr(h^{-1}vh^{-1}w),\\
\overline{tr(h^{-1}v)}  &  =h^{\bar{i}k}v_{k\bar{i}}
=tr(h^{-1}v).
\end{align*}
Therefore $\left\langle\cdot,\cdot\right\rangle_h^{\alpha}$ is real.

To show that $\left\langle\cdot,\cdot\right\rangle_h^{\alpha}$ is positive
definite, we need an important observation that will be used constantly in
this paper, that is, though $h^{-1}v$ is in general not a Hermitian matrix, it
can be canonically written as

\begin{equation}
h^{-1}v=h^{-\frac12}(h^{-\frac12}vh^{-\frac12})h^{\frac12},
\label{observation}
\end{equation}
where $h^{\frac12}\in Herm^+(r)$ is a matrix satisfying ($h^{\frac12}$)$^2=h$ and $h^{-\frac12}=(h^{\frac12})^{-1}$. In
\eqref{observation}, $u=h^{-\frac12}vh^{-\frac12}$ is a Hermitian
matrix, and it is positive definite whenever $v$ is.

When $\alpha=0$, we have $\left\langle v,v\right\rangle_h=tr(h^{-1}
vh^{-1}v)=tr(u^2)=\left\Vert u\right\Vert _F^2\geq0$, and $tr(u^2
)=0\Leftrightarrow u=0\Leftrightarrow v=0$.

In the general case (i.e. $\alpha>-1/r$), applying Jensen's inequality to the
function $x^2$ shows that the following inequality holds for every
diagonalizable real $r\times r$ matrix:

\[
tr(A)^2/r\leq tr(A^2),
\]
with equality obtained only when $A=cI$ for some $c\in
\mathbb{R}$. It follows that:

\begin{align*}
\left\langle v,v\right\rangle_h^{\alpha}  &  =tr(h^{-1}vh^{-1}v)+\alpha
tr(h^{-1}v)tr(h^{-1}v)\\
&  \geq(1/r+\alpha)tr(h^{-1}v)tr(h^{-1}v)\geq0
\end{align*}
and $\left\langle v,v\right\rangle_h^{\alpha}=0$ if and only if
$tr(h^{-1}v)=0$ and $h^{-1}v=cI$ for some $c\in
\mathbb{R}$, in other words, $v=0$.

Next, we move to $Herm^+(V)\subseteq Herm(V)$ for a fixed $r$-dimensional
complex vector space $V$. $Herm(V)$ (resp. $Herm^+(V)$) is the space of
Hermitian sesquilinear forms (resp. Hermitian inner products) on $V$. Clearly
we have $Herm(V)\cong Herm(r)$, $Herm^+(V)\cong Herm^+(r)$ and
$Herm^+(V)$ is an open convex cone in $Herm(V)$.

Before we proceed, we want to address some important isomorphisms that will
come into use later. For a fixed $h\in Herm^+(V)$, there is a conjugate $
\mathbb{C}$-linear isomorphism

\begin{equation}
\begin{aligned}
Sesqui(V)  &  \cong End(V)\\
v  &  \mapsto h^{-1}v
\end{aligned}
\label{iso1}
\end{equation}
Here $h^{-1}v\in End(V)$ is the one defined by the formula:

\[
v(\alpha,\beta)=h(\alpha,(h^{-1}v)\beta)\,\text{, for }\alpha,\beta\in V
\]

It is obvious that $v$ is Hermitian if and only if $h^{-1}v$ is self-adjoint
with respect to the inner product $h$, therefore the above isomorphism
restricts to an $
\mathbb{R}$-linear isomorphism:

\begin{equation}
\begin{aligned}
Herm(V)  &  \cong End(V,h)\\
v  &  \mapsto h^{-1}v
\end{aligned}
\label{iso2}
\end{equation}
where $End(V,h)$ is the real subspace of $End(V)$ formed by endomorphisms
that are self-adjoint w.r.t. $h$.

It helps to take a closer look at the isomorphism \eqref{iso2}. For $v\in
Herm(V)$, since $h^{-1}v$ is self-adjoint w.r.t. the inner product $h$, all
eigenvalues of $h^{-1}v$ are real. Moreover, when $v$ is positive definite
(resp. positive semidefinite), all eigenvalues of $h^{-1}v$ are positive
(resp. nonnegative) and therefore the determinant of $h^{-1}v$ is positive
(resp. nonnegative). This is proved by making use of the definition of $h^{-1}v$: if
$\lambda$ is an eigenvalue of $h^{-1}v$ with a nonzero eigenvector $\alpha$,
then $\lambda h(\alpha,\alpha)=h(\alpha,(h^{-1}v)\alpha)=v(\alpha,\alpha)>0$
(resp. $\geq0$).

In the bundle case, if we fix an $h\in\Gamma(Herm^+(E))$, then applying
\eqref{iso2} pointwisely gives an isomorphism of $\mathbb{R}$-vector bundles (here $End(E,h)_x=End(E_x,h_x)$):

\begin{equation}
Herm(E)\cong End(E,h)
\label{iso3}
\end{equation}
and an identification:

\begin{equation}
\begin{aligned}
\Gamma(Herm(E))  &  \cong\Gamma(End(E,h))\\
v  &  \mapsto h^{-1}v
\label{iso4}
\end{aligned}
\end{equation}

By the above discussion, at any point $x\in M$, the eigenvalues of
$h(x)^{-1}v(x)$ are all real. In this way, we obtain the eigenvalue functions
$\lambda_1\leq\cdots\leq\lambda_r$ associated with $h^{-1}v$, which are
smooth functions on $M$. When $v$ actually belongs to $\Gamma(Herm^+(E))$, we have $0<\lambda_1\leq\cdots\leq\lambda_r$; and when $v$ is merely
measurable or continuous, so is each $\lambda_i$.

Now back to our discussion of $Herm^+(V)$, for which we fix $h\in
Herm^+(V)$ and $v,w\in Herm(V)$. If we choice a basis $\varepsilon_1
,\dots,\varepsilon_r$ for $V$ and write $H_{i\bar{j}}=h(\varepsilon
_i,\varepsilon_{j})$, $V_{i\bar{j}}=v(\varepsilon_i,\varepsilon_{j}
)$, $W_{i\bar{j}}=w(\varepsilon_i,\varepsilon_{j})$, then by a simple calculation:

\[
v(\alpha,\beta)=a_iV_{i\bar{j}}\overline{b_{j}}=a_iH_{i\bar{l}
}(H^{\bar{l}k}V_{k\bar{j}}\overline{b_{j}})\text{, for
}\alpha=\varepsilon_ia_i,\beta=\varepsilon_{j}b_{j}
\]
we see that $(h^{-1}v)\beta=$ ($\overline{H^{\bar{l}k}V_{k\bar{j}}
}b_{j})\varepsilon_{l}$, thus the matrix representation of $h^{-1}v$ is
actually $\overline{(H^{\bar{l}k}V_{k\bar{j}})}$, i.e. the
conjugation of the matrix $H^{-1}V$ (recall that \eqref{iso1} is conjugate $\mathbb{C}$-linear).

The Riemannian metric on $Herm^+(V)$ is then defined as:

\begin{equation}
\begin{aligned}
\left\langle v,w\right\rangle_h^{\alpha}  &  \triangleq tr(h^{-1}v\cdot
h^{-1}w)+\alpha tr(h^{-1}v)tr(h^{-1}w)\\
&  =\overline{H^{\bar{k}i}H^{\bar{j}l}V_{i\bar{j}
}W_{l\bar{k}}}+\alpha\overline{H^{\bar{k}i}V_{i\bar{k}}H^{\bar{l}j}W_{j\bar{l}}}\\
&  =H^{\bar{k}i}H^{\bar{j}l}V_{i\bar{j}}
W_{l\bar{k}}+\alpha H^{\bar{k}i}V_{i\bar{k}
}H^{\bar{l}j}W_{j\bar{l}}.
\label{pointwisemetric}
\end{aligned}
\end{equation}

The subtle point is that, although the matrix expression of $h^{-1}v$ is not
$H^{-1}V$, but rather its conjugation, \eqref{pointwisemetric} still reads
$tr(H^{-1}VH^{-1}W)+\alpha tr(H^{-1}V)tr(H^{-1}V)$, for the conjugations are
all canceled out. As a result, the canonical isomorphism $Herm^+(V)\cong
Herm^+(r)$ w.r.t. the basis $\varepsilon_1,\dots,\varepsilon_r$ is an isometry.

Finally, we arrive at the bundle case. As mentioned in the Introduction, the
infinite-dimensional real vector space $\Gamma(Herm(E))$ is a Fr\'echet
space and the space of Hermitian metrics $\Gamma(Herm^+(E))$ is an open
convex cone in $\Gamma(Herm(E))$, making $\Gamma(Herm^+(E))$ trivially a
Fr\'echet manifold with its tangent space at each point canonically
identified with $\Gamma(Herm(E))$. In order to gain more generality, in the
bundle case, we assume that $\alpha>-1/r$ is not just a constant, but a smooth
function on $M$.

The natural Riemannian metric on $\Gamma(Herm^+(E))$ is defined by
integrating the scalar product \eqref{pointwisemetric}. More specifically, given
$h\in\Gamma(Herm^+(E))$ and $v,w\in\Gamma(Herm(E))\cong T_h\Gamma
(Herm^+(E))$, we set

\begin{equation}
(v,w)_h^{\alpha}\triangleq\int_Mtr(h^{-1}vh^{-1}w)+\alpha tr(h^{-1}v)tr(h^{-1}
w)\mu_g. 
\label{globalmetric}
\end{equation}
The function been integrated is $\left\langle v(x),w(x)\right\rangle
_{h(x)}^{\alpha(x)}=tr(h(x)^{-1}v(x)\cdot h(x)^{-1}w(x))+\alpha(x)tr(h(x)^{-1}
v(x))tr(h(x)^{-1}w(x))$, therefore its local expression is $h^{\bar{k}i}h^{\bar{j}l}v_{i\bar{j}}w_{l\bar{k}}+\alpha
h^{\bar{k}i}v_{i\bar{k}}h^{\bar{l}j}w_{j\bar{l}}$ (with respect to a smooth local frame). From the discussion in the
pointwise case, $(\cdot,\cdot)_h^{\alpha}$ is indeed a real inner product on
$\Gamma(Herm(E))$ and we denote the induced norm by $\left\Vert\cdot
\right\Vert _h^{\alpha}$, that is $\left\Vert\cdot\right\Vert _h^{\alpha
}=(v,v)_h^{\alpha}$ for $v\in\Gamma(Herm(E))$. Formally speaking, what we
are doing is giving a bundle Riemannian metric $\left\langle\cdot
,\cdot\right\rangle_h^{\alpha}$ on $Herm(E)$ and $(\cdot,\cdot)_h
^{\alpha}$ is the corresponding inner product on its space of smooth sections.

Throughout the paper, unless otherwise specified, we use the notation $d$ for
the distance function on $\Gamma(Herm^+(E))$ induced from the metric
\eqref{globalmetric} by taking the infimum of the lengths of paths between two given points.

The metric \eqref{globalmetric} is a weak Riemannian metric, which means that
the topology induced by $(\cdot,\cdot)_h^{\alpha}$ on $\Gamma(Herm(E))$ is
weaker than the smooth topology. This leads to some phenomena that are
unfamiliar from the world of finite-dimensional Riemannian geometry, or even
strong Riemannian metrics on Hilbert manifolds (\cite{L99}, \cite{Kl95}). For
instance, a weak Riemannian metric does not always induce a metric space
structure, in other words, $d$ is \textit{a priori} only a pseudometric, or even
vanishing (\cite{MM05}). Besides, for a weak Riemannian manifold, the metric
spray and the Levi-Civita connection may not exist (see, for example, Chap 4
in \cite{S23}). We shall see later in this paper the $L^2$ metrics as we
have defined them do not suffer from these pathologies.

There is a right action $(\phi^{\ast}h)(u,v)=h(\phi(u),\phi(v))$ ($\phi\in
GL(E),h\in\Gamma(Herm(E))$), which leaves $\Gamma(Herm^+(E))$ invariant.
Here $GL(E)$ is the group of vector bundle automorphisms of $E$ (covering the
identity map on $M$).  We can easily see that our metric is invariant under $GL(E)$, i.e. each $\phi\in GL(E)$ is an isometry:

\begin{equation}
\begin{aligned}
(\phi^{\ast}v,\phi^{\ast}w)_{\phi^{\ast}h}^{\alpha}  &  =\int_M\left\langle
(\phi^{\ast}v)(x),(\phi^{\ast}w)(x)\right\rangle_{(\phi^{\ast}h)(x)}
^{\alpha(x)}\mu_g(x)\\
&  =\int_M\left\langle v(x),w(x)\right\rangle_{h(x)}^{\alpha(x)}\,\mu
_g(x)=(v,w)_h^{\alpha}.
\label{acting}
\end{aligned}
\end{equation}
(This simple calculation has been omitted from \cite{K87}, so we include it here for completeness.)

\begin{remark}[About the ``complex gauge group'' $GL(E)$]
As mentioned above, $GL(E)$ is the group of vector bundle automorphisms of $E$; sometimes it is also denoted by $\mathscr{G}^\mathbb{C}$.
If the bundle is carried with a fixed Hermitian metric $h$, then the corresponding group should be the subgroup formed by automorphisms that preserve the Hermitian metric, and such a subgroup is usually denoted by $U(E,h)$, or simply by $\mathscr{G}$. These two groups played a role in various moduli spaces, such as the moduli space of holomorphic bundle structures, the moduli space of Hermitian-Yang-Mills connections, etc. The readers can consult \cite{K87}, Chap 7; \cite{D85} p.5; \cite{DK90}, Chap 6 and \cite{J91} p.54 for more information on this subject.

It is obvious that $GL(E)$ can be identified with the gauge transformation group of
the frame bundle of $E$ (which is a principal bundle). We want to remark that any natural Lie group
structure on $GL(E)$, such as the one introduced in \cite{ACMM89} or
\cite{W07}, necessarily makes the natural map $GL(E)\times E\rightarrow E$
smooth, thus our actions $GL(E)\times\Gamma(Herm(E))\rightarrow\Gamma
(Herm(E))$, $GL(E)\times\Gamma(Herm^+(E))\rightarrow\Gamma(Herm^+(E))$ are all smooth.
\end{remark}

\begin{remark}
There is a related flat structure on the space of Hermitian metric which we
now explain. Fixing an $h_0\in\Gamma(Herm^+(E))$ and a smooth parameter
function $\alpha>-1/r$, we have a fiber Riemannian metric $\left\langle
\cdot,\cdot\right\rangle_{h_0}^{\alpha}$ on $Herm(E)$ and, furthermore, an inner
product $(\cdot,\cdot)_{h_0}^{\alpha}$ on $\Gamma(Herm(E))$. Specifically,
it reads:

\begin{equation}
(v,w)_{h_0}^{\alpha}=\int_Mtr(h_0^{-1}vh_0^{-1}w)+\alpha tr(h_0
^{-1}v)tr(h_0^{-1}w)\mu_g. 
\label{flatmetric}
\end{equation}

We use $\left\Vert\cdot\right\Vert_{h_0}^{\alpha}$ to denote the induced
norm on $\Gamma(Herm(E))$. It is well-known that $(\Gamma(Herm(E)),\left\Vert
\cdot\right\Vert_{h_0}^{\alpha})$ is not complete, and its Hilbert
completion is the space of measurable sections of $Herm(E)$ subject to the
condition $\left\Vert u\right\Vert_{h_0}^{\alpha}<\infty$. The flat
Riemannian metric on $\Gamma(Herm^+(E))$ is defined by requiring that the
inner product on each tangent space $T_h\Gamma(Herm^+(E))\cong
\Gamma(Herm(E))$ be the same, namely \eqref{flatmetric}. It is clear that the
induced distance function is just $\left\Vert h_1-h_2\right\Vert_{h_0}
$ and shortest geodesics are straight line segments (note that $\Gamma
(Herm^+(E))$ is a convex domain in $\Gamma(Herm(E))$). Hence, endowed with
this flat Riemannian metric, $\Gamma(Herm^+(E))$ is not metrically complete,
and its metric completion is just its closure in the above mentioned Hilbert
completion of $\Gamma(Herm(E))$.
\end{remark}

\begin{remark}
When we talk about the gradient flow of a functional on a space, the space itself must be carried with some kind of a metric (\cite{V09} \S{23}). This is exactly the role of the natural metric \eqref{metric} in the study of Hermitian-Einstein structures. More specifically, the Hermitian-Yang-Mills flow on $\Gamma
(Herm^+(E))$ is $h_t^{-1}\dot{h_t}=-(i\Lambda F_t-cI)$. Donaldson's insight is that this parabolic flow is the gradient flow of an energy functional (now dubbed ``Donaldson's functional'') on $\Gamma(Herm^+(E))$, with respect to the metric \eqref{metric}. The readers can find the above discussion in \cite{K87} \S{6.3} or \cite{Siu87} p.40, although in the later book, the natural metric \eqref{metric} is not explicitly written down but implicitly used.

In \cite{Siu87} p.27, the gradient flow of another functional on $\Gamma(Herm^+(E))$, namely $J(h)=\int_M \vert i\Lambda F_h-cI\vert^2 \mu_g$, is shown to be $h_t^{-1}\dot{h_t}=-\Delta i\Lambda F_t$; in this case \eqref{metric} is again the background metric.
\end{remark}

\subsection{\texorpdfstring{$L^2$}{L2} metrics on the space of smooth sections of a fiber bundle}

The metric on $\Gamma(Herm^+(E))$ defined in the last subsection is a
special case of the so-called $L^2$ metric introduced in the Appendix of
\cite{FG89}. We now recall the definition. In this and the next subsection, $E$ is assumed
to be a smooth \textbf{fiber bundle} over the closed $n$-manifold $M$.

The setting is the following: we fix a smooth Riemannian metric on the
vertical tangent bundle $VTE$. This amounts to giving a smooth Riemannian
metric $\left\langle\cdot,\cdot\right\rangle^{x}$ on each fiber $E_x$ such
that $\left\langle\cdot,\cdot\right\rangle^{x}$ varies smoothly with respect
to $x\in M$. Also, in order to do the integration, we fix a Radon measure
$\mu$ on the base manifold $M$ (which could be induced by a smooth volume
form). All mapping spaces are endowed with the smooth compact-open topology.

It is well-known that the space of smooth sections $C^{\infty}(M,E)$ is a
smooth Fr\'echet manifolds with the tangent space at $\sigma\in$ $C^{\infty
}(M,E)$ identified with $\Gamma(\sigma^{\ast}(VTE))$ (cf. \cite{H82} p.85,
Example 4.1.2, we will work out the details later in the Appendix). For
$u,v\in\Gamma(\sigma^{\ast}(VTE))$, we set

\begin{equation}
\left(u,v\right)_{\sigma}\triangleq\int_M\left\langle u(x),v(x)\right\rangle_{\sigma(x)}^{x}\,\mu(x), 
\label{L2metric}
\end{equation}
where in this formula, $u(x)$ and $v(x)$ belong to $(\sigma^{\ast}
(VTE))_x=(VTE)_{\sigma(x)}=T_{\sigma\left(  x\right)  }E_x$ and
$\left\langle\cdot,\cdot\right\rangle_{\sigma(x)}^{x}$ is the inner product
on $T_{\sigma\left(x\right)  }E_x$ determined by the Riemannian metric
$\left\langle\cdot,\cdot\right\rangle^{x}$. More precisely, $\sigma$ gives a
fiber metric $\left\langle\cdot,\cdot\right\rangle_{\sigma}$ on the vector
bundle $\sigma^{\ast}(VTE)$, which further induces an $L^2$ inner product
$\left(\cdot,\cdot\right)_{\sigma}$ on $\Gamma(\sigma^{\ast}(VTE))$.

The formula \eqref{L2metric} defines a weak Riemannian metric on $C^{\infty
}(M,E)$, for the topology determined by $\left(\cdot,\cdot\right)_{\sigma
}$ on $\Gamma(\sigma^{\ast}(VTE))$ is the $L^2$ topology, not the smooth
topology. Again, the induced distance function is \textit{a priori} only a pseudometric.
In the next subsection we shall show that the distance function is indeed nondegenerate.

When there is a smooth real vector bundle $\xi$ over $M$ such that $E$ is an
open subbundle of $\xi$ in the category of smooth fiber bundles (cf.
\cite{P68}, Definition 12.2), the situation can be simplified. In this case,
$C^{\infty}(M,E)$ is open in the Fr\'echet space $\Gamma(\xi)$ (\cite{KM97}
p.300, Corollary 30.10 or \cite{H82} p.74, Example 3.17) so that each tangent
space of $C^{\infty}(M,E)$ is canonically identified with $\Gamma(\xi)$. Also,
$E_x$ is open in $\xi_x$, therefore the same thing happens for
$E_x\subseteq\xi_x$. Under this circumstance, in the formula
\eqref{L2metric}, $u,v$ are elements of $\Gamma(\xi)$, and $\left\langle
\cdot,\cdot\right\rangle_{\sigma(x)}^{x}$ is an inner product on $\xi_x$.
Note that since $(\sigma^{\ast}(VTE))_x=(VTE)_{\sigma(x)}=T_{\sigma\left(
x\right)  }E_x\cong\xi_x$, we have $\sigma^{\ast}(VTE)\cong\xi$.

The metric \eqref{globalmetric} on $\Gamma(Herm^+(E))$ fits perfectly well
in the above special situation, where each fiber $Herm^+(E_x)$ is endowed
with the Riemannian metric \eqref{pointwisemetric} with the parameter
$\alpha(x)$. Another example is the Ebin's metric on the space of Riemannian
metric (cf. the Appendix in \cite{FG89} and also, \cite{Cla13b}). A slightly
more complicated example is the weak Riemannian structure on the
diffeomorphism group. (For more information on the last one, the reader may
consult \cite{Sm07}, Section 5.)

\subsection{The distance function is nondegenerate}

In this subsection, we give a short account of the fact that the distance
function $d$ for the $L^2$ metric \eqref{L2metric} on $C^{\infty}(M,E)$ is
nondegenerate. As mentioned before, $d$ is \textit{a priori} only a pseudometric, so we
need to show that $d(\sigma_0,\sigma_1)\neq0$ for $\sigma_0\neq
\sigma_1$. This is done by proving an inequality (see below). The idea of
our proof is essentially contained in \cite{Cla10}, Proposition 17, but our
result seems to be simpler and more general.

We first note that there is a function $\Theta:$ $C^{\infty}(M,E)$ $\times
C^{\infty}(M,E)$ $\rightarrow
\mathbb{R}$ defined by

\[
\Theta(\sigma_0,\sigma_1)=\int_M\theta_x\left(\sigma_0(x),\sigma
_1(x)\right)\,\mu(x),
\]
where $\theta_x$ is the Riemannian distance function on the finite
dimensional manifold $E_x$. Clearly $\Theta$ is a metric on $C^{\infty
}(M,E)$ (in the sense of metric spaces). We have:

\begin{theorem}
$d(\sigma_0,\sigma_1)\geq\frac{1}{\sqrt{Vol(M,\mu)}}\Theta(\sigma
_0,\sigma_1)$.

\begin{proof}
Computing $d$ for arbitrary points $\sigma_0$ and $\sigma_1\in C^{\infty
}(M,E)$ involves an infinite-dimensional problem, since we have to find the
infimum of the expression:

\[
L(\sigma_{t})=\int_0^1\left\Vert\dot{\sigma}_{t}\right\Vert_{\sigma_{t}
}\,\mathrm{d}t=\int_0^1\left(\int_M\left\langle \dot{\sigma}_{t},\dot{\sigma}
_{t}\right\rangle_{\sigma_{t}}d\mu\right)^{1/2}\,\mathrm{d}t
\]
over all piecewise $C^1$-paths $\sigma_{t}$ connecting $\sigma_0$ and
$\sigma_1$. The key point here is to change the order of integration. From
H\"older's inequality,

\[
\int_M\sqrt{\left\langle\dot{\sigma}_{t},\dot{\sigma}_{t}\right\rangle
_{\sigma_{t}}}\,\mathrm{d}\mu\leq\left(\int_M\left\langle \dot{\sigma}_{t},\dot{\sigma}
_{t}\right\rangle_{\sigma_{t}}\,\mathrm{d}\mu\right)^{1/2}\sqrt{Vol(M,\mu)}.
\]

In other words,

\[
\left\Vert\dot{\sigma}_{t}\right\Vert_{\sigma_{t}}\geq\frac{1}
{\sqrt{Vol(M,\mu)}}\int_M\sqrt{\left\langle \dot{\sigma}_{t},\dot{\sigma
}_{t}\right\rangle_{\sigma_{t}}}\,\mathrm{d}\mu.
\]

It follows that

\begin{align*}
L(\sigma_{t})  &  =\int_0^1\left\Vert\dot{\sigma}_{t}\right\Vert
_{\sigma_{t}}\,\mathrm{d}t\geq\frac{1}{\sqrt{Vol(M,\mu)}}\int_0^1\int_M
\sqrt{\left\langle\dot{\sigma}_{t},\dot{\sigma}_{t}\right\rangle _{\sigma
_{t}}}\,\mathrm{d}\mu\,\mathrm{d}t\\
&  =\frac{1}{\sqrt{Vol(M,\mu)}}\int_M\left(\int_0^1\sqrt{\left\langle
\dot{\sigma}_{t}(x),\dot{\sigma}_{t}(x)\right\rangle_{\sigma_{t}(x)}^{x}
}\,\mathrm{d}t\right)\mathrm{d}\mu(x)\\
&  \geq\frac{1}{\sqrt{Vol(M,\mu)}}\int_M\theta_x(\sigma_0(x),\sigma
_1(x))\,\mathrm{d}\mu(x)=\frac{1}{\sqrt{Vol(M,g)}}\Theta(\sigma_0,\sigma_1),
\end{align*}
where the last ``$\geq$'' comes from the observation that for each $x\in M$,
$\sigma_{t}(x)(0\leq t\leq1)$ is a path in $E_x$ joining $\sigma_0(x)$ and
$\sigma_1(x)$.
\end{proof}
\end{theorem}

\begin{remark}
In \cite{Cla13b}, Theorem 2.1, it is proved that:

\[
d(\sigma_0,\sigma_1)\geq\left[\int_M\theta_x\left(\sigma_0(x),\sigma
_1(x)\right)^2\,\mu(x)\right]^{\frac12}.
\]

This result, together with H\"older's inequality, immediately gives the
above theorem. Nevertheless, we write down our proof since we find the idea
behind it is interesting.
\end{remark}

\subsection{Singular Hermitian metrics}

The concept of singular Hermitian metrics on a holomorphic line bundle was
first introduced by Demailly in \cite{DJP92} and has a key role in complex
geometry (see also, \cite{MM07} \S{2.3}). Singular Hermitian metrics on vector bundles were also introduced
and investigated in many papers, and by now its definition is pretty much
standard. We recall it in this subsection.

\begin{definition}
(\cite{BP08}, Section 3; \cite{R15}, Definition 1.1; \cite{PT18}, Definition
2.2.1; \cite{I20}, Definition 2.1) A singular Hermitian metric $h$ on $E$ is a
measurable map from the base manifold $M$ to the space of nonnegative
Hermitian forms on the fibers satisfying $0<\det h<\infty$ almost everywhere.
\end{definition}

This definition may seem a little vague. We think the one given in
\cite{HPS18} is clearer, so let's quote it here. It begins with the concept of
singular Hermitian inner products.

\begin{definition}
(\cite{HPS18}, Definition 16.1) A singular Hermitian inner product on a
finite-dimensional complex vector space $V$ is a function $\left\Vert\cdot\right\Vert$ : $V$ $\rightarrow\lbrack0,\infty]$ with the following properties:

(1) $\left\Vert\lambda v\right\Vert=\left\vert\lambda\right\vert\left\Vert v\right\Vert$.

(2) $\left\Vert v+w\right\Vert\leq\left\Vert v\right\Vert +\left\Vert
w\right\Vert$.

(3) $\left\Vert v+w\right\Vert^2+\left\Vert v-w\right\Vert^2
=2(\left\Vert v\right\Vert^2+\left\Vert w\right\Vert^2)$ (the
parallelogram law).
\end{definition}

The difference between a singular Hermitian inner product and a regular
one is, there may be nonzero vectors whose length is $0$,
and others whose length is $\infty$.

For a singular Hermitian inner product $\left\Vert\cdot\right\Vert$, both
$V_0=\{v\in V:\left\Vert v\right\Vert=0\}$ and $V_{\infty}=\{v\in
V:\left\Vert v\right\Vert<\infty\}$ are linear subspaces of $V$. We say
$\left\Vert\cdot\right\Vert$ is finite (resp. positive definite) if
$V_{\infty}=V$ (resp. $V_0=0$). Clearly $\left\Vert\cdot\right\Vert$ is
finite and positive definite if and only if $\left\Vert\cdot\right\Vert$ is
just a regular Hermitian inner product.

\begin{definition}
(\cite{HPS18}, Definition 17.1) A singular Hermitian metric $h$ on $E$ is an
assignment of a singular Hermitian inner product $h_x$ on each fiber $E_x$
in a measurable way. Also we assume that $h$ is finite and positive definite
almost everywhere.
\end{definition}

If we identify two singular Hermitian metrics that agree almost everywhere
on $M$, then the space of singular Hermitian metrics on $E$ is just the space
of measurable sections of $Herm^+(E)$, which we denote by $S(Herm^+(E))$.

For a smooth Hermitian metric $h$ on a holomorphic vector bundle $E$, the
concepts of Griffiths/Nakano seminegativity/semipositivity are defined using
the curvature tensor of the Chern connection (we shall not
recall these definitions). The only two things we need to know is that, in the
smooth case, $(E,h)$ is Griffiths seminegative if and only if$\ log\left\Vert u\right\Vert^2$ is plurisubharmonic for any local holomorphic section $u$;
and also, $(E,h)$ is Griffiths semipositive if and only if its dual $(E^{\ast
},h^{\ast})$ is Griffiths seminegative. These properties are crucial from our
point of view, since in the singular case, we can use these conditions as the
definition without referring to the curvature tensor.

\begin{definition}
(\cite{BP08}, Definition 3.1; \cite{R15}, Definition 1.2; \cite{PT18},
Definition 2.2.2; \cite{I20}, Definition 2.3) We say a singular Hermitian
metric $h$ (defined on a holomorphic vector bundle) is:

(1) Griffiths seminegative (or seminegatively curved in the sense of
Griffiths), if $\log\left\Vert u\right\Vert^2$, or equivalently, $\left\Vert
u\right\Vert^2$, is plurisubharmonic for any local holomorphic section $u$.

(2) Griffiths semipositive (or semipositively curved in the sense of
Griffiths), if the dual metric $h^{\ast}$ is a Griffiths seminegative
singular Hermitian metric on the dual bundle $E^{\ast}$.
\end{definition}

The definition of singular Hermitian metrics is far too liberal for us to
give the appropriate notion of ``curvature''.
Indeed, the curvature tensor of a singular Hermitian metric $h$ does not in
general make sense as a distribution, even when $h$ is Griffiths
seminegative. The example is found by Raufi in \cite{R15}, Theorem 1.5. We
will talk about this example later in Section 5.

By adding more conditions to the above definition, one can show the curvature
can be defined as a current with measure coefficients. This is the main result
of \cite{R15}.

\section{Basic geometry of $\Gamma(Herm^+(E))$}

\subsection{Metric spray, geodesics and curvature}

In this subsection we compute the various geometric objects associated with
our $L^2$ metric \eqref{globalmetric} on the space of Hermitian metrics. Our
starting point is the metric spray. The readers who are not familiar with the
concept ``metric spray'' should consult the
recent book \cite{S23}, Chap 4 or S. Lang's classic \cite{L99}. Roughly
speaking, for a weak Riemannian manifold, the metric spray may not exist; but
if it does exist, it is necessarily unique. The same thing holds for the
Levi-Civita connection. Also, the existence of the metric spray implies the
existence of the Levi-Civita connection, and in this case, they are related by
the formula given in \cite{S23} p.95, Proposition 4.24 or \cite{L99} p.199,
Theorem 2.1. In addition, geodesics and curvature can be computed directly
from the metric spray, so in some sense the metric spray is more convenient
than the Levi-Civita connection.

\begin{theorem}
The metric spray of the $L^2$ metric \eqref{globalmetric} exists and is
given by:

\[
B_h(v,w)=\frac12vh^{-1}w+\frac12wh^{-1}v,
\]
where $h$ takes its value in $\Gamma(Herm^+(E))$ and $v,w$ take their values
in $\Gamma(Herm(E))$. The expression $vh^{-1}w$ means that it is an element of
$\Gamma(Sesqui(E))$ satisfying $h^{-1}(vh^{-1}w)=h^{-1}vh^{-1}w\in
\Gamma(End(E))$ (here we are using the isomorphism \eqref{iso1}).

\begin{proof}
First, we must check that $\frac12vh^{-1}w+\frac12wh^{-1}v$ actually
belongs to $\Gamma(Herm(E))$. From the isomorphism \eqref{iso4}, we see that
$h^{-1}v$ and $h^{-1}w$ are both self-adjoint w.r.t. $h$, therefore
$h^{-1}vh^{-1}w+h^{-1}wh^{-1}v$ is also self-adjoint w.r.t. $h$. Again by
\eqref{iso4} we have $\frac12vh^{-1}w+\frac12wh^{-1}v\in
\Gamma(Herm(E))$.

The formula $(v,w)_h^{\alpha}=\int_Mtr(h^{-1}vh^{-1}w)+\alpha
tr(h^{-1}v)tr(h^{-1}w)\mu_g$ can be seen as a function with three variables
$v,w$ and $h$. Differentiate it in the $h$ direction using the definition of
the Bastiani calculus (\cite{H82} p.73, Definition 3.1.1; \cite{S23} p.7,
Definition 1.14), we obtain:

\begin{align*}
D(v,w)_h^{\alpha}(l)  &  =-\int_Mtr(h^{-1}lh^{-1}vh^{-1}w)\mu_g-\int
_Mtr(h^{-1}vh^{-1}lh^{-1}w)\mu_g\\
&  -\int_M\alpha tr(h^{-1}lh^{-1}v)tr(h^{-1}w)\mu_g-\int_M\alpha
tr(h^{-1}v)tr(h^{-1}lh^{-1}w)\mu_g.
\end{align*}

Then it follows that $\frac12D(v,v)_h^{\alpha}(l)-D(v,l)_h
^{\alpha}(v)=(vh^{-1}v,l)_h^{\alpha}$, which implies the symmetric bilinear
form $B_h(v,w)$ is the required metric spray.
\end{proof}
\end{theorem}

The geodesic equation can be deduced from the metric spray using the following formula
(cf. \cite{S23} p.92, Definition 4.19 and \cite{L99} p.101):

\begin{equation}
\ddot{\gamma}(t)=B_{\gamma(t)}(\dot{\gamma}(t),\dot{\gamma}(t)).
\label{spray-geodesic}
\end{equation}
In our case, it is then written as:

\[
\ddot{\gamma}(t)=\dot{\gamma}(t)\gamma(t)^{-1}\dot{\gamma}(t).
\]
By the convention we are using, this expression actually means that:

\[
\gamma(t)^{-1}\ddot{\gamma}(t)=\gamma(t)^{-1}\dot{\gamma}(t)\gamma(t)^{-1}
\dot{\gamma}(t),
\]
where both sides are elements of $\Gamma(End(E))$. In other words, the geodesic
equation in the space of Hermitian metrics reads:

\begin{equation}
\frac{d}{dt}\left(\gamma(t)^{-1}\dot{\gamma}(t)\right)=0.
\label{geodesicequation}
\end{equation}
This geodesic equation has been computed in \cite{K87} using the variation of energy.

The curvature tensor can also be explicitly written down via the metric spray
(cf. \cite{S23} p.100, Exercise 4.3.10 and \cite{L99} p.232, Proposition 1.2):

\begin{align*}
R_h(u,v)w  &  =B_h(B_h(v,w),u)-B_h(B_h(u,w),v)\\
&  +DB_h(u,w;v)-DB_h(v,w;u),
\end{align*}
here $u,v,w\in T_h\Gamma(Herm^+(E))\cong\Gamma(Herm(E))$ and $DB$ is the
derivative of $B$ in the $h$ direction.

A direct computation shows that $DB_h(v,w;l)=-\frac12vh^{-1}lh^{-1}w-\frac12wh^{-1}lh^{-1}v$. Then it follows that:

\begin{align*}
R_h(u,v)w  &  =-\frac14(wh^{-1}vh^{-1}u+uh^{-1}vh^{-1}w)\\
&  +\frac14(wh^{-1}uh^{-1}v+vh^{-1}uh^{-1}w).
\end{align*}

If we set $U=h^{-1}u$, $V=h^{-1}v$, $W=h^{-1}w\in\Gamma(End(E,h))$, then the
above formula can be written concisely as:

\begin{equation}
\begin{aligned}
h^{-1}R_h(u,v)w  &  =-\frac14(WVU+UVW-WUV-VUW)\\
&  =-\frac14\left[\left[U,V\right],W\right].
\label{curvature}
\end{aligned}
\end{equation}

Of course, the sectional curvature is easily derived from the curvature tensor
\eqref{curvature}. Given a pair of orthonormal vectors $u,v\in T_h
\Gamma(Herm^+(E))$, the sectional curvature of the plane spanned by them is

\begin{align*}
sec_h(u,v)  &  =(R_h(u,v)v,u)_h^{\alpha}\\
&  =\int_Mtr\left(h^{-1}R_h(u,v)vh^{-1}u\right)\mu_g\\
&  =\int_Mtr\left(\frac14\left[\left[U,V\right],V\right]U\right)\mu_g\\
&  =\frac14\int_Mtr\left(\left[U,V\right]^2\right)\mu_g\leq0.
\end{align*}
In the computation, the second term vanishes since $tr(h^{-1}R_h(u,v)v)=0$;
also, the ``$\leq$'' comes from the fact that
$\left[U,V\right]  $ is skew-Hermitian w.r.t. $h$.

In exactly the same way, we can study the metric \eqref{pointwisemetric} on
$Herm^+(V)$ and the metric \eqref{metric0} on $Herm^+(r)$. We just list
the results:

\begin{theorem}
For $Herm^+(V)$ $($resp. $Herm^+(r))$ , assuming that $h\in$ $Herm^+(V)$
$($resp. $Herm^+(r))$, $u,v,w\in Herm(V)$ $($resp. $Herm(r))$ $,U=h^{-1}u$,
$V=h^{-1}v$, $W=h^{-1}w$, we have:

metric spray :

\[
B_h(v,w)=\frac12vh^{-1}w+\frac12wh^{-1}v
\]

curvature tensor:

\[
h^{-1}R_h(u,v)w=-\frac14\left[\left[U,V\right],W\right]
\]

sectional curvature ($u,v$ are orthonormal):

\begin{equation}
sec_h(u,v)=\frac14tr\left(\left[U,V\right]^2\right)\leq0
\label{sec2}
\end{equation}

geodesic equation:

\begin{equation}
\frac{d}{dt}\left(\gamma(t)^{-1}\dot{\gamma}(t)\right)=0 
\label{geodesicequation2}
\end{equation}

The minor difference here is that in the $Herm^+(V)$ case, we are using the
isomorphism \eqref{iso2} implicitly so that $U,V,W$ and $\gamma(t)^{-1}
\dot{\gamma}(t)$ are elements of $End(V)$; while in the $Herm^+(r)$ case, the
operation $h^{-1}\cdot$ (or $\gamma(t)^{-1}\cdot$) is just the left multiplication of inverse matrix.
\end{theorem}

The consistency of the geometry of the bundle case and the fiberwise case is not a
coincidence, but rather the general feature of all the $L^2$ metrics on the
space of smooth sections (see the discussion in the Appendix of \cite{FG89}).
In \cite{FG89}, the proofs are not provided, so we shall say a little more
about this in our Appendix, especially the metric spray part. In view of this,
we actually only need to compute the fiberwise case then the bundle case will
just follow; but we feel that it is unnecessary and doesn't simplify the
computation, therefore we have shown the bundle case directly to make the
computation more concrete.

Next, we will look at the geodesics in $Herm^+(r)$ and $Herm^+(V)$ more
closely and write the distance function down explicitly. These results will be
used in the study of the metric completion of $\Gamma(Herm^+(E))$.

First, we have:

\begin{proposition}
The geodesic in $Herm^+(r)$ starting at $H$ with the initial velocity $A\in
Herm(r)$ is given by $\gamma(t)=He^{tH^{-1}A}=e^{tAH^{-1}}H$ ($t\in
\mathbb{R}$). In particular, $Herm^+(r)$ is geodesically complete and hence a
Cartan-Hadamard manifold (see \eqref{sec2}).

\begin{proof}
We must check that $\gamma(t)$ is indeed a curve in $Herm^+(r)$ first. From
\eqref{observation} we have $He^{tH^{-1}A}=H^{\frac12}e^{tH^{-\frac12
}AH^{-\frac12}}H^{\frac12}$, and it is an element of $Herm^+(r)$
since the naive exponential map maps $Herm(r)$ into $Herm^+(r)$.\ A
computation shows that $\dot{\gamma}(t)=Ae^{tH^{-1}A}=e^{tAH^{-1}}A$ and
$\gamma(t)^{-1}\dot{\gamma}(t)=H^{-1}A$, therefore $\gamma$ satisfies the
geodesic equation \eqref{geodesicequation2}.
\end{proof}
\end{proposition}

The above result shows that the exponential map for $Herm^+(r)$ at $I$ is
just the naive exponential map $\exp:$ $Herm(r)\rightarrow Herm^+(r)$, and it
is a diffeomorphism by the Cartan-Hadamard theorem. Clearly, the radial
geodesic segment $\gamma(t)=e^{tA}(0\leq t\leq1)$ is the unique shortest
geodesic joining $I$ to $e^{A}$, therefore $d(I,e^{A})=L(\gamma)=\sqrt
{tr(A^2)+\alpha tr\left(  A\right)  ^2}=\sqrt{\sum\lambda_i{}^2
+\alpha(\sum\lambda_i)^2}$, where $\lambda_1,\dots,\lambda_r$ are the
eigenvalues of $A$.

In general, for any $P,Q\in Herm^+(r)$, we have :

\begin{proposition}
(1). $d(P,Q)=\sqrt{\sum(\log\lambda_i)^2+\alpha(\log\lambda_1\cdots
\lambda_r)^2}$. ($\lambda_1,\dots,\lambda_r$ are the eigenvalues of
$P^{-1}Q$, they are all positive.)

(2). The shortest geodesic between P and Q is unique and is given by
$\gamma(t)=Pe^{tP^{-1}A}$ $(0\leq t\leq1)$, where $A$ is the unique $r\times
r$ Hermitian matrix s.t. $e^{P^{-1}A}=P^{-1}Q$.

\begin{proof}
From \eqref{observation} we have $P^{-1}Q=P^{-\frac12}P^{-\frac12
}QP^{-\frac12}P^{\frac12}$ and $P^{-\frac12}QP^{-\frac12}$ is
positive definite. Then there exists a unitary matrix $U$ s.t. $U^{-1}
P^{-\frac12}QP^{-\frac12}U=diag\{\lambda_1,\dots\lambda_r\}$. It
follows that $d(P,Q)=d(I,P^{-\frac12}QP^{-\frac12})=d(I,U^{-1}
P^{-\frac12}QP^{-\frac12}U)=\sqrt{\sum(\log\lambda_i)^2
+\alpha(\log\lambda_1\cdots\lambda_r)^2}$.

Also from \eqref{observation}, $e^{P^{-1}A}=P^{-1}Q$ can be rewritten as
$e^{P^{-\frac12}AP^{-\frac12}}=P^{-\frac12}QP^{-\frac12}$,
thus we can conclude the existence and uniqueness of such an $A$. A simple computation shows that $\gamma(t)=Pe^{tP^{-1}A}$ $(0\leq t\leq1)$ has length
$d(P,Q)$, therefore it is the unique shortest geodesic joining $P$ to $Q$ (the
uniqueness follows from the last proposition and the uniqueness of $A$).
\end{proof}
\end{proposition}

In the $Herm^+(V)$ case, using the isometry $Herm^+(V)\cong Herm^+(r)$
(fixing a basis $\varepsilon_1,\dots,\varepsilon_r$ for $V$), we have:

\begin{proposition}
(1). The geodesic in $Herm^+(V)$ starting at $h$ with the initial velocity
$v\in Herm(V)$ is given by $\gamma(t)=he^{th^{-1}v}$ ($t\in
\mathbb{R}
$). In particular, $Herm^+(V)$ is geodesically complete and hence a
Cartan-Hadamard manifold.

(2). $d(h_1,h_2)=\sqrt{\sum(\log\lambda_i)^2+\alpha(\log\lambda
_1\cdots\lambda_r)^2}$. ($\lambda_1,\dots,\lambda_r$ are the
eigenvalues of $h_1^{-1}h_2\in End(V)$, they are all positive.)

(3) The shortest geodesic between $h_1$ and $h_2$ is unique and is given
by $\gamma(t)=h_1e^{th_1^{-1}v}$ $(0\leq t\leq1)$, where $v$ is the unique
element in $Herm(V)$ s.t. $e^{h_1^{-1}v}=h_1^{-1}h_2$.

The only difference here is that again we are using the isomorphism
\eqref{iso2} implicitly so that $\gamma(t)=he^{th^{-1}v}$ actually means that
$h^{-1}\gamma(t)=e^{th^{-1}v}$.

\begin{remark}
Since a smooth curve $\gamma$ in $\Gamma(Herm^+(E))$ is a geodesic if and
only if $\gamma_x$ is a geodesic in $Herm^+(E_x)$ for all $x\in M$, we
know that $\Gamma(Herm^+(E))$ is geodesically complete and the typical
geodesic looks like $\gamma(t)=he^{th^{-1}v}$ ($t\in\mathbb{R}$) with $h\in\Gamma(Herm^+(E))$ and $v\in\Gamma(Herm(E))$ (this formula
actually means that $h^{-1}\gamma(t)=e^{th^{-1}v}$, via the map \eqref{iso4}).
\end{remark}
\end{proposition}

\subsection{The exponential map}

In the second part of this section, we study the exponential map for
$\Gamma(Herm^+(E))$ and show that it is always a diffeomorphism. We
emphasize that this is the infinite-dimensional case and the definition of
``diffeomorphism'' comes from the Bastiani
calculus (\cite{H82} p.73, Definition 3.1.1, \cite{S23} p.7, Definition 1.14)
or from the ``convenient calculus'' (\cite{KM97}) (they are equivalent in our case).

First, for any element $h\in Herm^+(V)$, if the basis $\varepsilon_1
,\dots,\varepsilon_r$ is chosen to be orthonormal w.r.t. $h$, then in the
isomorphism $Herm^+(V)\cong Herm^+(r)$ $h$ is mapped to the identity
matrix. In this case $exp_h:$ $Herm(V)\rightarrow Herm^+(V)$, the
exponential map at the point $h$, commutes with the naive exponential map
$\exp:$ $Herm(r)\rightarrow Herm^+(r)$:

\begin{equation}
\begin{tabular}
{lll}
$Herm(V)$ & $\cong$ & $Herm(r)$\\
$\downarrow exp_h$ & $\circlearrowright$ & $\downarrow \exp$\\
$Herm^+(V)$ & $\cong$ & $Herm^+(r)$
\end{tabular}
\label{pointwiseexp}
\end{equation}

Now fix an $h\in\Gamma(Herm^+(E))$ and denote the exponential map for
$\Gamma(Herm^+(E))$ at the base point $h$ by $Exp_h$. From the fact that a
smooth curve $\gamma$ in $\Gamma(Herm^+(E))$ is a geodesic if and only if
$\gamma_x$ is a geodesic in $Herm^+(E_x)$ for all $x\in M$, we know that
$Exp_h$ is defined on the entire tangent space $\Gamma(Herm(E))$ and is given by:

\begin{equation}
Exp_h(v)(x)=exp_{h(x)}(v(x)),
\label{E-e}
\end{equation}
where the right hand side is the exponential map of the finite dimensional
space $Herm^+(E_x)$, as in \eqref{pointwiseexp}.

We need a collected map:

\begin{equation}
\coprod_{x\in M}\exp_{h(x)}:Herm(E)\rightarrow Herm^+(E). \label{collectedexp}
\end{equation}

\begin{claim}
\eqref{collectedexp} is a smooth map.

\begin{proof}
Let $e_1,\dots,e_r$ be a smooth local frame for E that is orthonormal
w.r.t. the Hermitian metric $h$ and defined on the open set $U$. In this way
$Herm(E)\mid_{U}\cong U\times Herm(r)$ and $Herm^+(E)\mid$ $_{U}\cong
U\times Herm^+(r)$. According to \eqref{pointwiseexp}, the local expression
of \eqref{collectedexp} is the map sending $\left(  p,A\right)  $ to $\left(
p,e^{A}\right)  $, which is smooth.
\end{proof}
\end{claim}

The relation \eqref{E-e} means that $Exp_h:\Gamma(Herm(E))\rightarrow
\Gamma(Herm^+(E))$ is just the induced map of \eqref{collectedexp} on the
space of smooth sections. By \cite{KM97} p.300, Corollary 30.10 or \cite{H82}
p.85, Example 3.6.6, $Exp_h$ is smooth; the same argument applies to its
inverse. To sum up, we have:

\begin{theorem}
$Exp_h:\Gamma(Herm(E))\rightarrow\Gamma(Herm^+(E))$ is a diffeomorphism.
\end{theorem}

\section{The metric completion of $(\Gamma(Herm^+(E)),d)$}

In the last section the space of Hermitian metrics is shown to be geodesically
complete. However, we shall see in this section that it is not metrically
complete by writing down its metric completion concretely. Before that, we
remind the readers that our metric \eqref{globalmetric} is carried with a
smooth parameter function $\alpha>-1/r$.

\subsection{The integral representation for \texorpdfstring{$d$}{d}}

Our method is based on the following formula:

\begin{theorem}
For any $h_1,h_2\in\Gamma(Herm^+(E))$, we have:

\begin{equation}
d(h_1,h_2)=\left[\int_Md_x(h_1(x),h_2(x))^2\mu_g(x)\right]^{\frac12},
\label{mainformula}
\end{equation}
where $d_x$ is the distance function on $Herm^+(E_x)$ (endowed with the
Riemannian metric \eqref{pointwisemetric} with the parameter $\alpha(x)$).

\begin{proof}
The ``$\geq$'' part has already been proved in
\cite{Cla13b}, Theorem 2.1 so we only need to deal with the ``$\leq$'' part.

In Proposition 3.5 it is revealed that for any $x\in M$, $\gamma
_x(t)=h_1(x)e^{th_1(x)^{-1}v(x)}$ $(0\leq t\leq1)$ is the unique
shortest geodesic in $Herm^+(E_x)$ joining $h_1(x)$ to $h_2(x)$, where
$v(x)$ is the unique element in $Herm(E_x)$ s.t. $e^{h_1(x)^{-1}
v(x)}=h_1(x)^{-1}h_2(x)$. From this explicit expression, we know that
$v\in\Gamma(Herm(E))$ and $\gamma(t)=h_1e^{th_1^{-1}v}$ $(0\leq t\leq1)$
is a geodesic in $\Gamma(Herm^+(E))$ joining $h_1$ to $h_2$ (cf. Remark
3.6). Of course, at this point we don't know that $\gamma$ realizes the
distance between $h_1$ and $h_2$.

The length of $\gamma$, by definition, is the following:

\begin{align*}
L(\gamma)  &  =\int_0^1\left\Vert\dot{\gamma}(t)\right\Vert_{\gamma
(t)}\,\mathrm{d}t\\
&  =\int_0^1\left[\int_Mtr((\gamma_x(t)^{-1}\dot{\gamma}_x(t))^2
)+\alpha(x)tr(\gamma_x(t)^{-1}\dot{\gamma}_x(t))^2\mu_g(x)\right]^{\frac
12}\,\mathrm{d}t.
\end{align*}

Since $\gamma_x(t)$ is a geodesic in $Herm^+(E_x)$, it satisfies the
geodesic equation \eqref{geodesicequation2}, thus in the above expression
$\gamma_x(t)^{-1}\dot{\gamma}_x(t)$ is actually independent of $t$. The
following rather strange computation is valid:

\begin{align*}
L(\gamma)  &  =\int_0^1\left\Vert\dot{\gamma}(t)\right\Vert_{\gamma
(t)}\,\mathrm{d}t\\
&  =\int_0^1\left[\int_Mtr((\gamma_x(t)^{-1}\dot{\gamma}_x(t))^2
)+\alpha(x)tr(\gamma_x(t)^{-1}\dot{\gamma}_x(t))^2\mu_g(x)\right]^{\frac
{1}{2}}\,\mathrm{d}t\\
&  =\left[\int_Mtr((\gamma_x(t)^{-1}\dot{\gamma}_x(t))^2)+\alpha
(x)tr(\gamma_x(t)^{-1}\dot{\gamma}_x(t))^2\mu_g(x)\right]^{\frac12}\\
&  =\left[\int_M\left(\sqrt{tr((\gamma_x(t)^{-1}\dot{\gamma}_x(t))^2
)+\alpha(x)tr(\gamma_x(t)^{-1}\dot{\gamma}_x(t))^2}\right)^2\,\mu
_g(x)\right]^{\frac12}\\
&  =\left[\int_M\left(\int_0^1\sqrt{tr((\gamma_x(t)^{-1}\dot{\gamma}_x
(t))^2)+\alpha(x)tr(\gamma_x(t)^{-1}\dot{\gamma}_x(t))^2}\,\mathrm{d}t\right)^2\mu_g(x)\right]^{\frac12}\\
&  =\left[\int_ML(\gamma_x)^2\mu_g(x)\right]^{\frac12}\\
&  =\left[\int_Md_x(h_1(x),h_2(x))^2\mu_g(x)\right]^{\frac12}\leq
d(h_1,h_2).
\end{align*}

It follows that $d(h_1,h_2)=L(\gamma)=\left[\int_Md_x(h_1(x),h_2
(x))^2\mu_g(x)\right]^{\frac12}$.
\end{proof}
\end{theorem}

\begin{remark}
The above proof actually shows that any two points in $\Gamma(Herm^+(E))$
can be joined by a unique shortest geodesic segment that realizes the
distance (the uniqueness follows from Remark 3.6), therefore $(\Gamma
(Herm^+(E)),d)$ is a uniquely geodesic metric space (cf. \cite{BH13} p.4,
Definition 1.3). However, for a general infinite-dimensional geodesically
complete (weak or strong) Riemannian manifold, the minimizing geodesic between
two points may not exist; see, for example, \cite{S23} p.103, Example 4.43 or
the Remark in \cite{L99} p.226.
\end{remark}

\begin{remark}
For the general case $C^{\infty}(M,E)$ ($E$ is a fiber bundle over $M$, as in
Section 2.2 and 2.3), it may also happen that $d(\sigma_1,\sigma_2
)=[\int_M\theta_x(\sigma_1(x),\sigma_1(x))^2\mu_g(x)]^{\frac1
2}$. If for any point $x\in M$ we can find a shortest geodesic $\gamma
_x(t)$ $(0\leq t\leq1)$ in $E_x$ joining $\sigma_1(x)$ and $\sigma
_2(x)$\ such that for any fixed time $t$, the map $\gamma(t):M\rightarrow
E$ is smooth (Could this condition be fulfilled when each fiber is complete?),
then a similar computation runs like this:

\begin{align*}
L(\gamma)  &  =\int_0^1\left\Vert\dot{\gamma}(t)\right\Vert_{\gamma
(t)}\,\mathrm{d}t\\
&  =\int_0^1\left[\int_M\left\langle\dot{\gamma}_x(t),\dot{\gamma}
_x(t)\right\rangle _{\gamma_x(t)}^{x}\mu_g(x)\right]^{\frac12}\,\mathrm{d}t\\
&  =\left[\int_M\left\langle \dot{\gamma}_x(t),\dot{\gamma}_x(t)\right\rangle
_{\gamma_x(t)}^{x}\mu_g(x)\right]^{\frac12}\left(\left\langle
\dot{\gamma}_x(t),\dot{\gamma}_x(t)\right\rangle_{\gamma_x(t)}
^{x}\text{ is independent of }t\right)\\
&  =\left[\int_M\left(\sqrt{\left\langle\dot{\gamma}_x(t),\dot{\gamma}
_x(t)\right\rangle_{\gamma_x(t)}^{x}}\,\right)^2\mu_g(x)\right]^{\frac12}\\
&  =\left[\int_M\left(\int_0^1\sqrt{\left\langle\dot{\gamma}_x(t),\dot{\gamma
}_x(t)\right\rangle _{\gamma_x(t)}^{x}}\,\mathrm{d}t\right)^2\mu_g(x)\right]^{\frac12
}\left(\left\langle\dot{\gamma}_x(t),\dot{\gamma}_x(t)\right\rangle
_{\gamma_x(t)}^{x}\text{ is independent of }t\right)\\
&  =\left[\int_ML(\gamma_x)^2\mu_g(x)\right]^{\frac12}\\
& =\left[\int_M\theta_x(\sigma_1(x),\sigma_2(x))^2\mu_g(x)\right]^{\frac12}\leq d(\sigma_1,\sigma_2).
\end{align*}

Therefore $d(\sigma_1,\sigma_2)=L(\gamma)=[\int_M\theta_x(\sigma
_1(x),\sigma_1(x))^2\mu_g(x)]^{\frac12}$.
\end{remark}

We end this subsection by showing that the distance between two conformally
equivalent Hermitian metrics can be easily computed (for simplicity,
$\alpha>-1/r$ is assumed to be a constant):

\begin{corollary}
For any $h\in\Gamma(Herm^+(E))$ and any smooth functions $f,g$ on $M$, we have:

\[
d(e^{f}h,e^{g}h)=\sqrt{r(1+\alpha r)}\left\Vert f-g\right\Vert_2.
\]

\begin{proof}
$(e^{f(x)}h(x))^{-1}(e^{g(x)}h(x))=e^{g(x)-f(x)}I$, then by Proposition
3.5 (in this case each $\lambda_i=e^{g(x)-f(x)}$), $d_x(e^{f(x)}
h(x),e^{g(x)}h(x))=\sqrt{r(1+\alpha r)}\left\vert f(x)-g(x)\right\vert $, and
$d(e^{f}h,e^{g}h)=\sqrt{r(1+\alpha r)}\left\Vert f-g\right\Vert_2$. Another
way to see this is that, according to Remark 4.2, $\gamma(t)=e^{f}
he^{t(g-f)}(0\leq t\leq1)$ is the shortest geodesic segment joining $e^{f}h$
to $e^{g}h$, and a direct computation shows that $L(\gamma)=\sqrt{r(1+\alpha
r)}\left\Vert f-g\right\Vert_2$.
\end{proof}
\end{corollary}

\subsection{The space of \texorpdfstring{$L^2$}{L2} integrable singular Hermitian metrics as the
completion}

Recall that the space of singular Hermitian metrics (which is required to be
almost everywhere finite and positive definite) is denoted by $S(Herm^+
(E))$, which can be viewed as the space of measurable sections of
$Herm^+(E)$. Again we identify two singular Hermitian metrics that agree
almost everywhere on $M$.

Fixing an $h_0\in\Gamma(Herm^+(E))$, we set $L^2(Herm^+(E))$ to be:
\[
L^2(Herm^+(E))\triangleq\left\lbrace\sigma\in S(Herm^+(E)):\left[\int_Md_x
(\sigma(x),h_0(x))^2\mu_g(x)\right]^{\frac12}<\infty\right\rbrace.
\]

The triangle inequality and the Minkowski inequality together give:

\begin{align*}
&  \left[\int_Md_x(\sigma(x),h_0(x))^2\mu_g(x)\right]^{\frac12}\\
&  \leq\left[\int_Md_x(\sigma(x),h_1(x))^2\mu_g(x)\right]^{\frac12
}+\left[\int_Md_x(h_0(x),h_1(x))^2\mu_g(x)\right]^{\frac12}.
\end{align*}

Thus the definition of $L^2(Herm^+(E))$ is independent of the
choice of $h_0$ (the base point can actually be any element in
$L^2(Herm^+(E))$) and also $d_2(\sigma_1,\sigma_2)=[\int_M
d_x(\sigma_1(x),\sigma_1(x))^2\mu_g(x)]^{\frac12}$ defines a
metric space structure on $L^2(Herm^+(E))$. Clearly we have $L^2
(Herm^+(E))\supseteq C^{0}(Herm^+(E))\supseteq\Gamma(Herm^+(E))$. By our
main result \eqref{mainformula} in the last subsection, $d_2=d$ when
restricting to $\Gamma(Herm^+(E))$, hence from now on we will just write $d$
instead of $d_2$.

The elements in $L^2(Herm^+(E))$ will be referred to as $L^2$ integrable
singular Hermitian metrics.

In order to show $(L^2(Herm^+(E)),d)$ is the required metric completion,
we need: (1) $\Gamma(Herm^+(E))$ is dense in $(L^2(Herm^+(E)),d)$. (2)
$(L^2(Herm^+(E)),d)$ is a complete metric space.

The proof of (1) is based on the approximation method given in \cite{Ca23},
Theorem~2: first we show that $\Gamma(Herm^+(E))$ is dense in $C^{0}
(Herm^+(E))$ and then that $C^{0}(Herm^+(E))$ is dense in $L^2
(Herm^+(E))$. These facts are proved using the Tietze's Extension Theorem in
the version of \cite{DJ51}, Theorem 4.1 and classical approximation procedures
with partitions of unity. The only difference between our case and
\cite{Ca23}, Theorem 2 (where the Ebin's metric on the space of Riemannian
metrics is discussed) is that, in \cite{Ca23}, each fiber is incomplete, so an
extra step is required, namely to construct a bundle with each fiber being the
completion of the original finite-dimensional fiber, then the completion of
the space of Riemannian metrics is the space of $L^2$ integrable sections of
this new bundle; while in our case, each fiber $Herm^+(E_x)$ is already
complete so things are much simpler. Since we do not need to construct a new
bundle, the argument in the first paragraph of \cite{Ca23}, Theorem 2 is not
needed and meanwhile the second paragraph of his proof can be translated
almost word by word to be the proof of (1) (we will not repeat it here). So
the only thing left to be proved is (2). We start with two lemmas.

\begin{lemma}
Let $\{h_{k}\}$ be a sequence in $L^2(Herm^+(E))$ s.t. $\sum_{k=1}^{\infty}d(h_{k},h_{k+1})<\infty$, then the function
$\Omega(x)=$ $\sum_{k=1}^{\infty}d_x\left(h_{k}(x),h_{k+1}
(x)\right)$ is $L^2$ integrable i.e. $\Omega\in L^2(M,\mu_g)$.

\begin{proof}
The Minkowski inequality gives
\begin{align*}
&  \left[\int_M \left(\sum_{k=1}^N d_x\left(h_{k}(x),h_{k+1}
(x)\right)\right)^2\mu_g(x)\right]^{\frac12}\\
&  \leq\sum_{k=1}^N\left[\int_M(d_x(h_{k}(x),h_{k+1}
(x))^2\mu_g(x)\right]^{\frac12}\\
&  \leq\sum_{k=1}^{\infty} d(h_{k},h_{k+1})<\infty.
\end{align*}
Then the conclusion follows from Beppo Levi's monotone convergence theorem.
\end{proof}
\end{lemma}

\begin{lemma}
Let $\{h_{k}\}$ be a sequence in $L^2(Herm^+(E))$ s.t. $\sum_{k=1}^{\infty} d(h_{k},h_{k+1})<\infty$, assume in addition that
$h_{k}$ converges to a singular Hermitian metric $\sigma$ a.e. on $M$.

Then $\sigma\in L^2(Herm^+(E))$ and $\{h_{k}\}$ converges to $\sigma$ in
$L^2(Herm^+(E))$, too.

\begin{proof}
Set $g_{k}(x)=d_x(h_{k}(x),\sigma(x))$ and we already know that
$g_{k}\rightarrow0$ a.e. on $M$. Note that $g_k(x)\leq\sum_{l=k}^{{N-1}} d_x(h_{l}(x),h_{l+1}(x))+d_x(h_n
(x),\sigma(x))$. As $N$ tends to $\infty$, the first term is controlled by
$\Omega$ and the second term goes to zero, therefore $g_{k}\leq\Omega\in
L^2(M,\mu_g)$. It follows that $g_{k}\in L^2(M,\mu_g)$ and
$g_{k}\rightarrow0$ in $L^2(M,\mu_g)$ (by the $L^2$-version of the
dominated convergence theorem), in other words, $\sigma\in L^2(Herm^+(E))$
and $h_{k}\rightarrow\sigma$ in $L^2(Herm^+(E))$.
\end{proof}
\end{lemma}

Finally, we have:

\begin{theorem}
$(L^2(Herm^+(E)),d)$ is complete.

\begin{proof}
Given a Cauchy sequence $\{h_{k}\}$ in $L^2(Herm^+(E))$, by passing to a
subsequence, we can assume that $\sum_{k=1}^{\infty}d(h_{k},h_{k+1})<\infty$. Define $\Omega$ as above. Since it is $L^2$
integrable, $\Omega$ is almost everywhere finite, which implies, for almost
every $x\in M$, $h_{k}(x)$ is a $d_x$-Cauchy sequence in $Herm^+(E_x)$.
We have seen in Section 3 that $Herm^+(E_x)$ is complete, therefore for
almost every $x\in M$, $h_{k}(x)$ converges to some element $\sigma(x)$ in
$Herm^+(E_x)$. In this way we get a singular Hermitian metric $\sigma$ and
the above lemma shows that $\{h_{k}\}$ converges to $\sigma$ in $L^2
(Herm^+(E))$.
\end{proof}
\end{theorem}

So far we have proved that $(L^2(Herm^+(E)),d)$ is the completion of
$\Gamma(Herm^+(E))$. Next, we show:

\begin{theorem}
$(\Gamma(Herm^+(E)),d)$ is a CAT(0) space. (Recall that it has already been proved
to be a uniquely geodesic metric space, see Remark 4.2.)

\begin{proof}
First, for any finite-dimensional simply-connected complete Riemannian
manifold, being nonpositively curved is equivalent to being CAT(0) (cf.
\cite{BH13} p.159), therefore each $Herm^+(E_x)$ is CAT(0). Given any
geodesic triangle $\Delta(g,h,k)$ in $\Gamma(Herm^+(E))$, from our knowledge
of shortest geodesic in $\Gamma(Herm^+(E))$ (Remark 4.2), we know that for
any $x\in M$, $\Delta(g,h,k)(x)$ is a geodesic triangle in $Herm^+(E_x)$.
Now set $\bar{g}(x)=(0,0)$, $\bar{h}(x)=(d_x(g(x),h(x)),0)$ and let $\bar
{k}\left(  x\right)  $ be the point in $\mathbb{R}^2$ s.t. $\left\vert\bar{k}\left(x\right)-\bar{g}(x)\right\vert
=d_x(k(x),g(x))$, $\left\vert\bar{k}\left(x\right)-\bar{h}
(x)\right\vert=d_x(k(x),h(x))$, then $\Delta(\bar{g}(x),\bar{h}(x),\bar
{k}\left(x\right))\subseteq\mathbb{R}
^2$ is a comparison triangle for $\Delta(g,h,k)(x)$.

$L^2(M,\mathbb{R}
^2)$, endowed with the norm $\left\Vert f\right\Vert_2=(\int
_M\left\vert f\right\vert^2\mu_g)^{\frac12}$, is a complete CAT(0)
space (\cite{BH13} p.167). Thus to show that $\Gamma(Herm^+(E))$ satisfies
the CAT(0) inequality, it suffices to take a comparison triangle in $L^2(M,\mathbb{R}
^2)$ rather than $\mathbb{R}
^2$. We claim that $\Delta(\bar{g},\bar{h},\bar{k})\subseteq L^2(M,\mathbb{R}
^2)$ is a comparison triangle for $\Delta(g,h,k)$. Indeed, using our main
result \eqref{mainformula}, we have $\left\Vert\bar{g}-\bar{h}\right\Vert
_2=(\int_M\left\vert\bar{g}-\bar{h}\right\vert^2\mu_g)^{\frac12
}=(\int_Md_x(g(x),h(x))^2\mu_g(x))^{\frac12}=d(g,h)$ and
similarly, $\left\Vert\bar{g}-\bar{k}\right\Vert_2=d(g,k)$, $\left\Vert
\bar{h}-\bar{k}\right\Vert_2=d(h,k)$.

Now for any points $p,q\in\Delta(g,h,k)$, we have $p(x),q(x)\in\Delta
(g,h,k)(x)$. Let $\bar{p}(x),\bar{q}(x)\in$ $\Delta(\bar{g}(x),\bar{h}
(x),\bar{k}\left(x\right))$ be the corresponding points in the comparison
triangle. The CAT(0) property of $Herm^+(E_x)$ yields $d_x
(p(x),q(x))\leq\left\vert\bar{p}(x)-\bar{q}(x)\right\vert$, then again by
\eqref{mainformula} we have $d(p,q)\leq\left\Vert\bar{p}-\bar{q}\right\Vert_2$. The conclusion follows.
\end{proof}
\end{theorem}

Since the metric completion of a CAT(0) space is also a CAT(0) space
(\cite{BH13} p.187, Corollary 3.11), we have:

\begin{corollary}
$(L^2(Herm^+(E)),d)$ is a CAT(0) space.
\end{corollary}
\begin{remark}
The computation in Corollary 4.4 can obviously be carried over to
$L^2(Herm^+(E))$ (we assume $\alpha>-1/r$ is a constant).

Firstly, given an $h\in L^2(Herm^+(E))$ and a measurable function $f$ on
$M$, we have $d(e^{f}h,h)=\sqrt{r(1+\alpha r)}\left\Vert f\right\Vert _2$,
thus in this case $e^{f}h$ is $L^2$ integrable if and only if $f\in
L^2(M,\mu_g)$. Secondly, given $f,g\in L^2(M,\mu_g)$, we have the same
formula $d(e^{f}h,e^{g}h)=\sqrt{r(1+\alpha r)}\left\Vert f-g\right\Vert_2$
as in Corollary 4.4.
\end{remark}

\subsection{An alternative description of the metric completion}

In this short subsection we closely follow \cite{Ca23}, Theorem 3 to give an
alternative description of the metric completion that might be useful.
However, this construction is only valid when $\alpha$ is a constant.

We use $\rho$ to denote the distance function on $Herm^+(r)$ induced by
the Riemannian metric \eqref{metric0} with the parameter $\alpha$. Of course, in view of the isometry $Herm^+(r)\cong Herm^+(E_x)$, $\rho$ is just the distance function $d_x$ introduced in Theorem 4.1.  

Fixing a $A\in Herm^+(r)$, define:

\[
L^2(M,Herm^+(r))\triangleq\left\lbrace f:M\rightarrow Herm^+(r):\left[\int_M
\rho(f(x),A)^2\mu_g(x)\right]^{\frac12}<\infty\right\rbrace,
\]
which is independent of the choice of $A$. From the triangle inequality and
the Minkowski inequality, we know for any $f,g\in L^2(M,Herm^+(r))$:

\begin{align*}
&  \left[\int_M\rho(f(x),g(x))^2\mu_g(x)\right]^{\frac12}\\
&  \leq\left[\int_M\rho(f(x),A)^2\mu_g(x)\right]^{\frac12}+\left[\int_M
\rho(g(x),A)^2\mu_g(x)\right]^{\frac12}<\infty
\end{align*}
and

\begin{align*}
&  \left[\int_M\rho(f(x),g(x))^2\mu_g(x)\right]^{\frac12}\\
&  \leq\left[\int_M\rho(f(x),h(x))^2\mu_g(x)\right]^{\frac12}+\left[\int
_M\rho(h(x),g(x))^2\mu_g(x)\right]^{\frac12}.
\end{align*}
Hence $d_2(f,g)=[\int_M\rho(f(x),g(x))^2\mu_g(x)]^{\frac12}$ is a
distance function on $L^2(M,Herm^+(r))$.

\begin{proposition}
$(L^2(Herm^+(E)),d)\cong(L^2(M,Herm^+(r)),d_2)$.

\begin{proof}
We fix a triangulation for $M$ and denote by $V_1,\cdots,V_{k}$ the
interiors of the maximal dimensional simplices. Also we can assume that each
$V_i$ is contained in some local trivialization chart for the bundle
$Herm^+(E)$, thanks to the barycentric subdivision of a simplicial complex
and the existence of the Lebesgue number for an open cover. Therefore we have
$\phi_i:Herm^+(E)\mid_{V_i}\cong V_i\times Herm^+(r)\rightarrow
Herm^+(r)$. For any $\sigma\in L^2(Herm^+(E))$, the $Herm^+(r)$-valued
function $f_{\sigma}(x)=\phi_i\circ\sigma(x)(x\in V_i)$ is a.e. defined on
$M$ and:
\begin{align*}
\left\lbrack\int_Md_x(\sigma(x),\tau(x))^2\mu_g(x)\right\rbrack^{\frac12}  & =\left\lbrack\sum_i\int_{V_i}d_x(\sigma(x),\tau(x))^2\mu_g(x)\right\rbrack^{\frac12}\\
& =\left\lbrack\sum_i\int_{V_i}\rho(f_{\sigma}(x),f_{\tau}(x))^2\mu_g(x)]\right\rbrack^{\frac12}\\
&  =\left\lbrack\int_M\rho(f_{\sigma}(x),f_{\tau}(x))^2\mu_g(x)\right\rbrack^{\frac12}.
\end{align*}

Hence $\sigma\mapsto f_{\sigma}$ provides an isometric isomorphism between
$(L^2(Herm^+(E)),d)$ and $(L^2(M,Herm^+(r)),d_2)$.
\end{proof}
\end{proposition}

\section{The holomorphic case}

In this section, $E$ is assumed to be a holomorphic vector bundle over the
compact complex manifold $M$, and still $M$ is carried with a fixed Riemannian
metric $g$.

\subsection{Relation to Griffiths seminegativity/semipositivity}

The goal of this subsection is to prove that, a Griffiths
seminegative/semipositive singular Hermitian metric is always $L^2$
integrable. First, we shall give some convenient equivalent
conditions of ``$L^2$ integrable''.

Fix an $h_0\in\Gamma(Herm^+(E))$ throughout out this subsection. Given any
singular Hermitian metric $h\in S(Herm^+(E))$, we write $H=h_0^{-1}h$,
which is a measurable section of the vector bundle $End(E,h_0)$. Since all
eigenvalues of $H$ are positive, we obtain the eigenvalue functions
$0<\lambda_1^H\leq\cdots\leq\lambda_r^H$ (they are measurable
functions defined on $M$ since $h$ is measurable). As mentioned before, each
$\lambda_i^H$ is continuous/smooth when $h$ is continuous/smooth.

By definition, $h\in L^2(Herm^+(E))$ if and only if $\int_M
d_x(h(x),h_0(x))^2\mu_g(x)<\infty$, and this is equivalent to the
condition that $\sum(\log\lambda_i^H)^2+\alpha(\log\det H)^2$ is
integrable on $M$, due to the computation given in Proposition 3.5. We can
take a step further:

\begin{proposition}
The following three conditions are equivalent for a singular Hermitian metric $h\in
S(Herm^+(E))$:

(1) $h\in L^2(Herm^+(E))$.

(2) $\log\lambda_1^H,\dots,\log\lambda_r^H\in L^2(M,\mu_g)$.

(3) $\log\lambda_1^H$ and $\log\lambda_r^H\in L^2(M,\mu_g)$.

And any of these implies:

(4) $\log\det H\in L^2(M,\mu_g)$.
\begin{proof}
(2) and (3) are obviously equivalent and they imply (4) since
$\log\det H=\log\lambda_1^H\cdots\lambda_r^H=\log\lambda_1^H
+\cdots+\log\lambda_r^H$. Also, it is clear that (2)\&(4)$\implies$(1) and
(1)\&(4)$\implies$(2), so the only thing left to be proved is that
(1)$\implies$(4). Applying Jensen's inequality to the function $x^2$ shows
that:
\[
\frac{(\log\det H)^2}{r^2}=\left(\frac{\log\lambda_1^H+\cdots+\log\lambda
_r^H}{r}\right)^2\leq\frac{(\log\lambda_1^H)^2+\cdots+(\log\lambda_r
^H)^2}{r},
\]
therefore $\sum(\log\lambda_i^H)^2+\alpha(\log\det H)^2\geq(\frac1
r+\alpha)(\log\det H)^2$, which gives the result.
\end{proof}
\end{proposition}

Next we need the concepts of ``good'' local
trivialization chart and ``bounded'' singular
Hermitian metrics.

\begin{definition}
For a local trivialization chart $(U,\phi)$ for $E$, if there is a larger
local trivialization chart $(V,\psi)$ such that $\bar{U}\subseteq V$ and
$\psi\mid$ $_{U}=\phi$, then we say $(U,\phi)$ is ``good''.
\end{definition}

The compactness of $M$ implies that $M$ can be covered by finitely many ``good'' local trivialization charts.

\begin{definition}
For a singular Hermitian metric $h$, we say $h$ is bounded, if, on any
``good'' local trivialization chart
$(U,\phi)$, the matrix expression of $h$ on $U$, say $(h_{i\bar{j}})$, is
bounded, i.e. there exists a constant $C>0$, s.t. $\left\vert h_{i\bar
{j}}\right\vert\leq C$ on $U$ for $1\leq i,j\leq r$.
\end{definition}

Clearly, any continuous Hermitian metric is bounded.

When $h$ is bounded, we have better results:

\begin{proposition}
For a bounded singular Hermitian metric $h$, the function $\lambda_r^H$ is
bounded from above by some constant.

\begin{proof}
Assume that $U$ is a ``good'' local
trivialization chart. Since $M$ can be covered by finitely many
``good'' local trivialization charts, we only
need to show $\lambda_r^H$ is bounded from above on $U$. Let us consider
the Frobenius norm of the local matrix expression of $h_0$, i.e. $\left\Vert
h_0(x)\right\Vert_{F}$. Thanks to the positive-definiteness and the
smoothness of $h_0$, $\left\Vert h_0(x)\right\Vert_{F}$ is bounded from
below by some positive constant $\delta$ on $U$. For any $x\in U$, there is a
nonzero vector $v$ in $E_x$ s.t. $h_0(x)^{-1}h(x)v=\lambda_r^H(x)v$,
then we have $\left\vert\lambda_r^H(x)\right\vert\left\Vert v\right\Vert
_2\leq\left\Vert h_0(x)^{-1}h(x)\right\Vert_{F}\left\Vert v\right\Vert
_2$ and hence $\left\vert\lambda_r^H(x)\right\vert\leq\left\Vert
h_0(x)^{-1}h(x)\right\Vert_{F}=\frac{\left\Vert h(x)\right\Vert_{F}
}{\left\Vert h_0(x)\right\Vert_{F}}\leq\frac{C}{\delta}$.
\end{proof}
\end{proposition}

\begin{proposition}
For a bounded singular Hermitian metric $h$, (4) implies (2) (hence also implies (1)).

\begin{proof}
The last proposition shows that $\lambda_r^H$ is bounded from above by some
constant $C$. Also, from $\det H=\lambda_1^H\cdots\lambda_r^H\leq
\lambda_1^H(\lambda_r^H)^{r-1}$ we know $\lambda_1^H\geq
\frac{\det H}{(\lambda_r^H)^{r-1}}\geq$ $\frac{\det H}{C^{r-1}}$. Therefore
each $\log\lambda_i^H$ is bounded from below by an $L^2$ integrable
function, provided that (4) is satisfied.
\end{proof}
\end{proposition}

\begin{remark}
So far we haven't say anything about the holomorphic condition yet, thus the
above discussions are all valid in the smooth category. However, the following
results require that $E$ and $M$ are holomorphic.
\end{remark}

\begin{theorem}
Any Griffiths seminegative singular Hermitian metric is $L^2$ integrable.

\begin{proof}
It is shown in Lemma 2.2.4 of \cite{PT18}, that a Griffiths seminegative
singular Hermitian metric is always bounded, so all we need to prove is that
$\log\det H\in L^2(M,\mu_g)$. In any local holomorphic trivialization chart,
$\log\det H=\log\det h-\log\det h_0$ ($h$ and $h_0$ are local matrix expressions
w.r.t. the holomorphic local frame), where $h_0$ is smooth so it is
irrelevant. From \cite{R15}, Proposition 1.3 or \cite{HPS18}, Proposition 25.1
we know $\log\det h$ is plurisubharmonic. On the other hand, in \cite{I20}, Lemma
2.6, it is proved that any plurisubharmonic function is locally square
integrable. The conclusion follows.
\end{proof}
\end{theorem}

If $h$ is a Griffiths semipositive singular Hermitian metric on $E$, then by
definition, $h^{\ast}$ is a Griffiths seminegative singular Hermitian metric
on $E^{\ast}$; also, $h_0^{\ast}$ is a smooth Hermitian metric on $E^{\ast}
$. Since the dual metric is locally represented by the transport of the
inverse matrix of the original one, we see the eigenvalues of $(h_0^{\ast
})^{-1}h^{\ast}$ are $0<\frac{1}{\lambda_r^H}\leq\cdots\leq\frac
{1}{\lambda_1^H}$. Therefore the above theorem implies:

\begin{corollary}
Any Griffiths semipositive singular Hermitian metric is $L^2$ integrable.
\end{corollary}

We end this subsection by providing two examples that may help to illustrate
our results:

\begin{example}
This is the example given in \cite{R15}, Theorem 1.5. Let $\Delta\subseteq
\mathbb{C}
$ be the unit disk, $E=\Delta\times
\mathbb{C}
^2$ be the trivial bundle of rank 2 over $\Delta$. Let $h$ be the singular
Hermitian metric on $E$ represented by the matrix $($using the global frame
$(1,0),(0,1))$
\[
\left(
\begin{array}{cc}
1+\left\vert z\right\vert^2 & z\\
\bar{z} & \left\vert z\right\vert^2
\end{array}
\right).
\]
As shown in Raufi's paper, $h$ is Griffiths seminegative, but the
``curvature'' $\Theta_h$ is not a current
with measure coefficients. We shall see that $h$ is $L^2$ integrable.
(Though the base manifold is not compact, it's a bounded domain endowed with
the Lebesgue measure, so our previous discussions are still valid.) Let
$h_0$ be the smooth Hermitian metric on $E$ represented by the identity
matrix, then the matrix expression of $H=h_0^{-1}h\in S(End(E))$ $($which is
just the space of measurable maps from $\Delta$ to $End(
\mathbb{C}
^2))$ is:
\[
\left(
\begin{array}{cc}
1+\left\vert z\right\vert^2 & \bar{z}\\
z & \left\vert z\right\vert^2
\end{array}
\right),
\]
where the conjugation comes from the discussion above \eqref{pointwisemetric}.

The eigenvalue functions of $H$ are $\lambda_1(z)=\lambda_2(z)=\left\vert
z\right\vert^2$, while $\left\vert\log\left\vert z\right\vert
^2\right\vert^2=4(\log\frac{1}{\left\vert z\right\vert})^2$, which is
integrable on $\Delta$. Therefore the condition (2) in Proposition 5.1 is satisfied.
\end{example}

\begin{example}
Let $L$ be a holomorphic line bundle over a compact complex manifold $M$
(endowed with a fixed Riemannian metric). If $h$ is a singular Hermitian
metric on $L$ and $s$ is a nonvanishing holomorphic section of $L$ defined on
an open set $U$, then with respect to $s$, the metric $h$ is represented by a
measurable function $e^{\phi}=\left\vert s\right\vert_h^2:U\rightarrow
\lbrack0,\infty]$. By definition, $h$ is Griffiths seminegative if and only if
all such $\phi$ are plurisubharmonic (which is consistent with the original
definition given in \cite{DJP92}).

We conclude that $h$ is $L^2$ integrable if and only if all such $\phi$ are
square-integrable. Indeed, a fixed smooth Hermitian metric $h_0$ is
represented by a positive smooth function $e^{\phi_0}$, thus the local
expression of $H=h_0^{-1}h$ is $e^{\phi-\phi_0}$. It follows from
Proposition 5.1 that $h$ is $L^2$ integrable if and only if all such
$\phi-\phi_0$ are square-integrable, but $\phi_0$ is irrelevant since it is smooth.

When $L=M\times\mathbb{C}$ is the trivial line bundle, the situation is the simplest:

(1) A Griffiths seminegative singular Hermitian metric on $L$ is just a
plurisubharmonic function on $M$.

(2) An $L^2$ integrable singular Hermitian metric on $L$ is just an $L^2$
integrable function on $M$, and in this case, $d(\phi_1,\phi_2
)=\sqrt{1+\alpha}\left\Vert\phi_1-\phi_2\right\Vert_2$ (cf. Remark
4.10., here we assume $\alpha>-1$ is a constant).
\end{example}

\subsection{The space of Hermitian metrics on a complex manifold}

In this subsection we briefly discuss the space of Hermitian metrics on a
compact complex manifold $M$. By definition, a Hermitian metric on $M$ is either

(a) a $J$-invariant Riemannian metric $g$ ($J$ is the almost complex
structure), or

(b) a Hermitian metric $h$ on the holomorphic tangent bundle $T^{1,0}M$.

These two definitions are equivalent. To get (b) from (a), we first
extend $g$ to be a $\mathbb{C}
$-\textbf{bilinear} symmetric form on the complexified tangent bundle $TM\otimes
\mathbb{C}
$ and again denote it by $g$. It is easy to see that: (1) $g(\bar{v}
,\bar{w})=$ $\overline{g(v,w)}$ (2) $g$ is still $J$-invariant ($J$ has also
been extended) (3) $g(v,w)=0$ when $v,w\in T^{1,0}M$ or $v,w\in T^{0,1}M$ (due
to (2)). In particular, $g$ is determined by its value on $T^{1,0}M\times
T^{0,1}M$. It follows that $h(v,w)=g(v,\bar{w})$ is the (b) we want. On the other hand, we can obtain (a) from (b) by going through the above procedure backward.

There is an isomorphism of complex vector bundles:

\begin{align*}
	(TM,J) &  \cong T^{1,0}M\\
	v &  \longmapsto\sigma(v)=\frac12(v-iJv)
\end{align*}
Via this map, we have $g(v,w)=2Reh(v,w)$.
\begin{remark}
In some literature, $g$ and $h$ are related by $g(v,w)=Reh(v,w)$. Our convention requires the constant 2.
\end{remark}
From the above discussion, the space of Hermitian metrics on $M$ is just the space $\Gamma(Herm^+(T^{1,0}M))$, and it is a subset of
$R(M)$ ($R(M)$ is the space of all Riemannian metrics on $M$). Indeed, we have $\Gamma(Herm(T^{1,0}M))\hookrightarrow\Gamma
(S^2T^{\ast}M)$ and it restricts to $\Gamma(Herm^+(T^{1,0}
M))\hookrightarrow R(M)$. 

Moreover, the following result is valid:

\begin{theorem}
The space of Hermitian metrics on $M$, $\Gamma(Herm^+(T^{1,0}M))$, is a split submanifold of $R(M)$.

\begin{proof}
For the definition of (split) submanifold in infinite dimension, the readers
may consult \cite{S23} p.17, Definition 1.35 or \cite{N06} p.313, Definition 1.3.5.
	
We set:
\[
\Gamma^{-1}(S^2T^{\ast}M,J)=\{l\in\Gamma(S^2T^{\ast}M):l(Jv,Jw)=-l(v,w)\}.
\]

Clearly, $\Gamma^{-1}(S^2T^{\ast}M,J)$ is a subspace of $\Gamma(S^2T^{\ast
}M)$. If we extend $l$ to be a $
\mathbb{C}
$-bilinear symmetric form on $TM\otimes
\mathbb{C}
$ (in the same way as above), we see $l(v,w)=0$ when $v\in T^{1,0}M,w\in T^{0,1}M$ or $v\in T^{0,1}M,w\in T^{1,0}M$, therefore $l$ is determined by its value on $T^{1,0}M\times
T^{1,0}M$, which gives us $\Gamma^{-1}(S^2T^{\ast}M,J)\cong\Gamma
(S^2(T^{1,0}M)^{\ast})$.
	
Every $g\in\Gamma(S^2T^{\ast}M)$ can be uniquely written as $g(u,v)=\frac
{g(u,v)+g(Ju,Jv)}{2}+\frac{g(u,v)-g(Ju,Jv)}{2}$, where the first term is
$J$-invariant and the second term vanishes only when $g$ itself is
$J$-invariant. It follows that $\Gamma(S^2T^{\ast}M)=\Gamma(Herm(T^{1,0}
M))\oplus\Gamma^{-1}(S^2T^{\ast}M,J)$, from which we know $\Gamma(Herm^+(T^{1,0}M))$ is a split submanifold of $R(M)$.
\end{proof}
\end{theorem}

It seems that the most natural weak Riemannian metric on $\Gamma
(Herm^+(T^{1,0}M))$ is given by:

\begin{equation}
(v,w)_h\triangleq\int_Mtr\left(h^{-1}vh^{-1}w\right)\mu_h, 
\label{2Ebin}
\end{equation}
where $h\in\Gamma(Herm^+(T^{1,0}M))$ and $v,w\in\Gamma(Herm(T^{1,0}M))$, and
$\mu_h$ is the volume element associated with $h$.

At this point we need a linear injection:

\begin{align*}
M(n,\mathbb{C})  &  \longrightarrow M(2n,\mathbb{R})\\
A  &  \longmapsto\left(
\begin{array}{cc}
	ReA & ImA\\
	-ImA & ReA
\end{array}
\right)  =\tau(A)=\hat{A}
\end{align*}
whose restriction is a Lie group homomorphism from $GL(n,\mathbb{C})$ to $GL(2n,\mathbb{R})$.

If $h,v,w$ are represented by the Hermitian matrix $H,V,W$ w.r.t. a
holomorphic local coordinate $z_1,\dots,z_n$, then their Riemannian
counterparts are represented by $2\hat{H},2\hat{V},2\hat{W}$, respectively, hence:

\begin{align*}
tr(\hat{H}^{-1}\hat{V}\hat{H}^{-1}\hat{W})  &  =tr(\tau(H^{-1}VH^{-1}W))\\
&  =2Retr(H^{-1}VH^{-1}W)\\
&  =2tr(H^{-1}VH^{-1}W).
\end{align*}

It follows that, up to a constant factor 2, \eqref{2Ebin} is just the
restriction of Ebin's metric to $\Gamma(Herm^+(T^{1,0}M))$. (Even if we use
the convention given in Remark 5.1, the factor 2 still pops up.) Moreover,
all computations for the Ebin's metric on $R(M)$ can be directly carried over
to $\Gamma(Herm^+(T^{1,0}M))$. In particular, $\Gamma(Herm^+(T^{1,0}M))$
is a totally geodesic submanifold of $R(M)$.

\subsection{Sobolev metrics on the space of Hermitian metrics}

One may go a little further to consider the Sobolev metrics on the space of
Hermitian metrics. For a fixed $h\in\Gamma(Herm^+(E))$, we have already
defined a bundle Riemannian metric $\left\langle\cdot,\cdot\right\rangle_h$ on $Herm(E)$ (see the paragraph following \eqref{globalmetric}), which,
combined with $g$,\ gives each $\otimes^kT^{\ast}M\otimes Herm(E)$ a
Riemannian structure $\left\langle\cdot,\cdot\right\rangle_{h,g}$.

On the other hand, the Chern connection of $h$ naturally induces a connection
$\nabla_h$ on $Herm(E)$, which further extends to an antiderivation (i.e.
$\nabla_h(\alpha\wedge\beta)=d(\alpha)\wedge\beta+(-1)^{deg(\alpha)}
\alpha\wedge\nabla_h(\beta)$):

\[
\nabla_h:\Omega^k(M,Herm(E))\rightarrow\Omega^{k+1}(M,Herm(E)).
\]

Now for any positive integer $s$, there is an inner product on $\Gamma
(Herm(E))$:

\begin{equation}
(u,v)_{h,g}^s\triangleq\sum_{j=0}^s\int_M\left\langle \nabla_h^ju,\nabla_h^jv\right\rangle_{h,g}
\mu_g. 
\label{Hsmetric}
\end{equation}
When $s=0$, it is just the $L^2$ metric \eqref{metric}.

\eqref{Hsmetric} is still a weak Riemannian metric on $\Gamma(Herm^+(E))$,
for the topology determined by \eqref{Hsmetric} on the tangent space
$\Gamma(Herm(E))$ is the $C^s$ topology, not the smooth topology; another
way to see this is that, the Hilbert completion of $\Gamma(Herm(E))$ w.r.t.
the inner product \eqref{Hsmetric} is the Sobolev $H^s$ space of sections,
customarily written as $H^s(Herm(E))$.

\begin{remark}
There is a second way to define $(u,v)_{h,g}^s$, namely, the connection
$\nabla_h$ on $Herm(E)$ and the Levi-Civita connection $\nabla_g$ together
form a tensored connection $\nabla=\nabla_g\otimes\nabla_h$ on each
$\otimes^kT^{\ast}M\otimes Herm(E)$, and we can use $\nabla$ instead of
$\nabla_h$ in the definition of \eqref{Hsmetric}.
\end{remark}

The family of weak Riemannian metrics $(\cdot,\cdot)_{h,g}^s(s\in
\mathbb{N}
)$ is a so-called graded Riemannian structure on the Fr\'echet manifold
$\Gamma(Herm^+(E))$ (cf. \cite{DR19}, Definition 3.19). We hope to
investigate these metrics in a future paper, especially when $s=1$.

\section{Appendix: The exponential map of the \texorpdfstring{$L^2$}{L2} metric on the space of
smooth sections}

We have seen at the end of Section 3, that the exponential map of the metric
\eqref{globalmetric} is a diffeomorphism. The same thing holds for the
exponential map of Ebin's metric on the space of Riemannian metrics
(\cite{GMM91}, Theorem 3.4). This is not a coincidence but rather a
consequence of the nonpositive curvature. In this Appendix we investigate the
exponential map of the general $L^2$ metric introduced in Section 2.2.
Specifically, we use the Nash-Moser theorem to prove that, if each fiber is
nonpositively curved, then the exponential map is a local diffeomorphism. This
result can be seen as a version of Cartan-Hadamard theorem in the tame
Fr\'echet setting. Before we get started, remember that the space of the
smooth sections of a vector bundle (resp. fiber bundle) is a tame Fr\'echet
space (resp. tame Fr\'echet manifold) (\cite{H82} p.139, Corollary 1.3.9 and
p.146, Theorem 2.3.1).

We shall first treat the special case, that is, the fiber bundle $E$ is an
open subbundle of the real vector bundle $\xi$ (see Section 2.2), for in this
case there is a global coordinate for $C^{\infty}(M,E)$ so that the proofs are
more straightforward. After that, we move on to the general case and see
what must be altered.

As before, we start with the metric spray. A direct computation yields:

\begin{proposition}
The metric spray for the $L^2$ metric on $C^{\infty}(M,E)$ exists and is
given by:
\begin{equation}
B(\sigma;u,v)(x)=B^{x}(\sigma(x);u(x),v(x)),x\in M
\label{spray-special}
\end{equation}
where $\sigma$ takes its value in $C^{\infty}(M,E)$, $u,v$ take their values
in $\Gamma(\xi)$, and $B^{x}:E_x\times\xi_x\times\xi_x\rightarrow\xi
_x$ is the metric spray for the finite-dimensional Riemannian manifold
$E_x$.
\end{proposition}

Using this proposition and the formula \eqref{spray-geodesic}, we see that a
smooth curve $\gamma$ in $C^{\infty}(M,E)$ is a geodesic if and only if
$\gamma_x$ is a geodesic in $E_x$ for all $x\in M$, which is a statement
given in \cite{FG89} Corollary A.4. Other results in the Appendix of
\cite{FG89} can also be deduced from the metric spray \eqref{spray-special}.

For a fixed $\sigma\in C^{\infty}(M,E)$, the exponential map at $\sigma$ is
denoted by $Exp_{\sigma}:\Omega\rightarrow C^{\infty}(M,E)$, where
$\Omega\subseteq$ $\Gamma(\xi)$ is the maximal domain of its definition. On
the other hand, let $exp_{\sigma(x)}:\varepsilon_x\subseteq\xi
_x\rightarrow E_x$ be the exponential map of the finite-dimensional
Riemannian manifold $E_x$ with the base point $\sigma(x)$. Collecting
$exp_{\sigma(x)}$ together, we obtain a map:
\begin{equation}
\coprod_{x\in M} \exp_{\sigma(x)}:\varepsilon=\coprod_{x\in M}\varepsilon_x\longrightarrow\coprod_{x\in M}E_x=E 
\label{collected}
\end{equation}

Since $\sigma$ is a smooth object and the family of Riemannian metrics
$\left\langle\cdot,\cdot\right\rangle^{x}$ varies smoothly w.r.t. $x\in M$,
$\varepsilon$ is open in $\xi$ and \eqref{collected} is a smooth map. In fact,
the metrics $\left\langle\cdot,\cdot\right\rangle^{x}$ (i.e. the given
metric on the vertical tangent bundle $VTE$) give rise to a unique bundle spray
over $E$ and \eqref{collected} is just the corresponding exponential map. (For
the concept of the ``bundle spray'', see
\cite{P68},\ Chap 12.)

Our understanding of the geodesics in $C^{\infty}(M,E)$ shows that
$\Omega=\{v\in\Gamma(\xi):v(x)\in\varepsilon_x$ for all $x\in M\}$ and
$Exp_{\sigma}(v)(x)=exp_{\sigma(x)}(v(x))$, i.e. $Exp_{\sigma}$ is the induced
map of \eqref{collected} on the space of smooth sections. By \cite{KM97}
p.300, Corollary 30.10 or \cite{H82} p.74, Example 3.1.7, $\Omega$ is open in
$\Gamma(\xi)$ and $Exp_{\sigma}:\Omega\rightarrow C^{\infty}(M,E)$ is smooth;
and by \cite{H82} p.145, Theorem 2.2.6, $Exp_{\sigma}$ is tame. In summary,
we have:

\begin{proposition}
$Exp_{\sigma}=(\coprod_{x\in M}\exp_{\sigma(x)})_{\ast}:\Omega\rightarrow C^{\infty}(M,E)$ is smooth and tame.
\end{proposition}

At this point, we can stop for a second to look at a particularly special case,
that is, each $E_x$ is a Cartan-Hadamard manifold. In this case,
$\varepsilon=\xi$ and \eqref{collected} is a diffeomorphism. Applying
\cite{KM97} p.300, Corollary 30.10 or \cite{H82} p.74, Example 3.1.7 to both
$Exp_{\sigma}$ and its inverse, we get:

\begin{proposition}
If each $E_x$ is a Cartan-Hadamard manifold, then $Exp_{\sigma}:\Gamma
(\xi)\rightarrow C^{\infty}(M,E)$ is a diffeomorphism.
\end{proposition}

Our metric \eqref{globalmetric} on the space of Hermitian metrics is an
example of such a very special case. However, the Ebin's metric on the space
of Riemannian metrics is not, since each fiber is nonpositively curved but
incomplete. Nevertheless, it falls into the realm of the
next result:

\begin{theorem}
If each $E_x$ is nonpositively curved, then $Exp_{\sigma}:\Omega\rightarrow
C^{\infty}(M,E)$ is a local diffeomorphism.

\begin{proof}
Since each $E_x$ is nonpositively curved, the differential of its
exponential map (at the point it is defined) is a linear isomorphism. (This
result does not require the geodesic completeness; see the proof of Theorem
3.7 in \cite{L99} p.251-252.) Therefore the differential $D(exp_{\sigma
(x)}):\varepsilon_x\times\xi_x\rightarrow\xi_x$ has the property that
$D(exp_{\sigma(x)})(v,\cdot):\xi_x\rightarrow\xi_x$ is a linear
isomorphism \textbf{for all} $v\in\varepsilon_x$. We write $G_x
(v,w)=D(exp_{\sigma(x)})(v,\cdot)^{-1}(w):\varepsilon_x\times\xi
_x\rightarrow\xi_x$ (which is also a smooth map).

Now we consider $D(Exp_{\sigma}):\Omega\times\Gamma(\xi)\rightarrow\Gamma
(\xi)$. By definition,

\begin{equation}
\begin{aligned}
D(Exp_{\sigma})(f,v)(x)  &  =\lim_{t\to0}\frac{Exp_{\sigma
}(f+tv)-Exp_{\sigma}(f)}{t}(x)\\
&  =\lim_{t\to0}\frac{exp_{\sigma(x)}(f(x)+v(x))-exp_{\sigma
(x)}(f(x))}{t}\\
&  =D(exp_{\sigma(x)})(f(x),v(x)).
\end{aligned}
\label{compu-special}
\end{equation}

In other words, $D(Exp_{\sigma})=\left(\coprod_{x\in M}
D(exp_{\sigma(x)})\right)_{\ast}$, hence is smooth and tame.

From this expression we know that, \textbf{for any} $g\in\Omega$,
$D(Exp_{\sigma})(g,\cdot):\Gamma(\xi)\rightarrow\Gamma(\xi)$ is invariable
and the family of inverses $V\left(Exp_{\sigma}\right):\Omega\times
\Gamma(\xi)\rightarrow\Gamma(\xi)$ is given by $\left(\coprod_{x\in M} G_x\right)_{\ast}$, which is also smooth and tame.

The conclusion then follows from the Nash-Moser inverse function theorem, as
stated in \cite{H82} p.171-172.
\end{proof}
\end{theorem}
\begin{remark}
The key point in the proof is that, the nonpositive curvature condition makes
$D(Exp_{\sigma})(g,\cdot)$ invertible \textbf{for all} $g\in\Omega$, not
just one; and $Exp_{\sigma}$, $D(Exp_{\sigma})$ and the whole collection of
inverses $V\left(Exp_{\sigma}\right)$ are all smooth and tame for being
the case considered in \cite{H82} p.74, Example 3.1.7 and p.145, Theorem
2.2.6. As a consequence, $Exp_{\sigma}$ satisfies all the requirements needed
in the Nash-Moser inverse function theorem.
\end{remark}

Next we move to the general case, where $E$ is just assumed to be a smooth
fiber bundle over $M$ (without being an open subbundle of a vector bundle). In
this case, $C^{\infty}(M,E)$ is a Fr\'echet manifold, not just an open
subset of a Fr\'echet space, so there is no global coordinate on $C^{\infty
}(M,E)$. In order to deal with this Fr\'echet manifold case, we need to know
the explicit coordinate chart around any given point $f\in C^{\infty}(M,E)$. In
\cite{H82}, Hamilton just wrote that, as we quote, ``it is
easy to construct a diffeomorphism from a neighborhood of the zero section of
$f^{\ast}(VTE)$ to a neighborhood of the image of $f$ in $E$''
, but we feel that this construction is not so obvious (at least for the
author), so here we provide the details.

The diffeomorphism we need is actually constructed in the proof (not the
result) of Theorem 12.10, ``existence theorem for
VBN'', in \cite{P68}. The metric on $VTE$ gives rise to a
metric spray on each fiber $E_x$. Collecting them together we get a bundle spray over
$E$ (so the bundle spray is not randomly chosen, but determined by the metric
on $VTE$). By \cite{P68}, Lemma 12.9, there are strictly positive smooth
functions $\lambda$ and $\mu$ on $E$, such that $exp_{e}$ (the exponential map
for $e\in E_x$) maps $B(\lambda(e))\subseteq T_{e}E_x$
diffeomorphically onto its image, which contains $B(e,\mu(e))\subseteq E_x$.
The set
\[
\Lambda\triangleq\left\lbrace(x,v)\in f^{\ast}(VTE):\left\vert v\right\vert <\lambda(f(x))\right\rbrace
\]
is open in $f^{\ast}(VTE)$ and $F=\coprod_{x\in M} exp_{f(x)}:\Lambda\rightarrow E$ is fiber-preserving (on each fiber, it is
just $exp_{f(x)}:B(\lambda(f(x)))\rightarrow E_x$). Then $U=Im(F)$ is an
open set in $E$ containing $f(M)$. As a result:

\begin{proposition}
$F$ is the required diffeomorphism and $(\widetilde{U},(F^{-1})_{\ast})$
serves as a coordinate chart around $f$, where $\widetilde{U}=\{s\in
C^{\infty}(M,E):s(M)\subseteq U\}$ (following the notation used in \cite{H82}
p.74, Example 3.1.7).
\end{proposition}

This explains why the tangent space at $f$ is identified with $\Gamma(f^{\ast
}(VTE))$.

Using the above coordinate, we can describe the metric spray:

\begin{proposition}
The metric spray for the $L^2$ metric on $C^{\infty}(M,E)$ exists and its
local expression is given by:
\[
B(g;u,v)(x)=G_{f}^{-1}\left[B^{x}(g(x);G_{f}(u(x)),G_{f}(v(x)))\right],x\in M
\]
where $g$ takes its value in $\widetilde{U}$, $u,v$ take their values in
$\Gamma(f^{\ast}(VTE))$, $B^{x}$ is the metric spray for $E_x$, and
$G_{f}=d(exp_{f(x)})_{exp_{f(x)}^{-1}(g(x))}$ is the map from $f^{\ast
}(VTE)_x$ to $T_{g(x)}E_x\cong g^{\ast}(VTE)_x$.
\end{proposition}

All results in the Appendix of \cite{FG89} can be derived from the above proposition.

We then turn our attention to the exponential map in this general case. For a
fixed $\sigma\in C^{\infty}(M,E)$, the exponential map for $C^{\infty}(M,E)$
(resp. $E_x$) is still written as $Exp_{\sigma}:\Omega\rightarrow C^{\infty
}(M,E)$ (resp. $exp_{\sigma(x)}:\varepsilon_x\rightarrow E_x$), but this
time $\Omega\subseteq$ $\Gamma(\sigma^{\ast}(VTE))$ (resp. $\varepsilon
_x\subseteq T_{\sigma\left(  x\right)  }E_x\cong(\sigma^{\ast}(VTE))_x
$). The collected smooth map still looks like \eqref{collected}, although in
this case, $\varepsilon$ is an open subset of $\sigma^{\ast}(VTE)$.

Again we have $\Omega=\{v\in\Gamma(\sigma^{\ast}(VTE)):v(x)\in\varepsilon_x$
for all $x\in M\}$ and $Exp_{\sigma}(v)(x)=exp_{\sigma(x)}(v(x))$, i.e.
$Exp_{\sigma}=\left(\coprod_{x\in M}exp_{\sigma(x)}\right)_{\ast}$. But we cannot conclude that $Exp_{\sigma}$ is
smooth and tame at this point, for it is not the vector bundle case. However,
this is not really a problem in view of the above discussion. Fix an
$l\in\Omega$ and let $f=Exp_{\sigma}(l)$. Let $W$ be the preimage of $U$ under
the map $\coprod_{x\in M} exp_{\sigma(x)}$, then we obtain:
\begin{equation}
W\longrightarrow U \longrightarrow \Lambda
\label{cop1}
\end{equation}

The first arrow is the map $\coprod_{x\in M} exp_{f(x)})^{-1}$. The induced map on the smooth sections is:
\begin{equation}
\widetilde{W}\xrightarrow{Exp_{\sigma}}\widetilde{U}\xrightarrow{(F^{-1})_{\ast}}\widetilde{\Lambda}
\label{cop2}
\end{equation}

Therefore the general case is reduced to the vector bundle case by using the
explicit local coordinate $(\widetilde{U},(F^{-1})_{\ast})$. As before,
\cite{KM97} p.300, Corollary 30.10 or \cite{H82} p.74, Example 3.1.7 \& p.145,
Theorem 2.2.6 implies:

\begin{proposition}
(1) $\Omega$ is open in $\Gamma(\sigma^{\ast}(VTE))$ and $Exp_{\sigma}
:\Omega\rightarrow C^{\infty}(M,E)$ is smooth and tame.

(2) If each $E_x$ is a Cartan-Hadamard manifold, then $Exp_{\sigma}
:\Gamma(\sigma^{\ast}(VTE))\rightarrow C^{\infty}(M,E)$ is a diffeomorphism.
\end{proposition}

The only thing left to be proved is the following theorem, and in the process
the coordinate $(\widetilde{U},(F^{-1})_{\ast})$ is vital, too.

\begin{theorem}
If each $E_x$ is nonpositively curved, then $Exp_{\sigma}:\Omega\rightarrow
C^{\infty}(M,E)$ is a local diffeomorphism.

\begin{proof}
Since each $E_x$ is nonpositively curved, the differential of its
exponential map (at the point it is defined) is a linear isomorphism.
Following the notation given above, we choose an arbitrary $l\in\Omega$, set
$f=Exp_{\sigma}(l)$, and focus our attention on the subset $W$ and
$\widetilde{W}$. The differential $D(exp_{f(x)}^{-1}\circ exp_{\sigma
(x)}):W_x\times\sigma^{\ast}(VTE)_x\rightarrow f^{\ast}(VTE)_x$ has the
property that $D(exp_{f(x)}^{-1}\circ exp_{\sigma(x)})(v,\cdot):\sigma^{\ast
}(VTE)_x\rightarrow f^{\ast}(VTE)_x$ is a linear isomorphism \textbf{for
all} $v\in W_x$. We write $G_x(v,w)=D(exp_{f(x)}^{-1}\circ exp_{\sigma
(x)})(v,\cdot)^{-1}(w):W_x\times f^{\ast}(VTE)_x\rightarrow\sigma^{\ast
}(VTE)_x$.

Now consider $D((F^{-1})_{\ast}\circ Exp_{\sigma}):\widetilde{W}\times
\Gamma(\sigma^{\ast}(VTE))\rightarrow\Gamma(f^{\ast}(VTE))$ (the differential
of \eqref{cop2}). A computation similar to \eqref{compu-special} yields
$D((F^{-1})_{\ast}\circ Exp_{\sigma})=\left(\coprod_{x\in M}D(exp_{f(x)}^{-1}\circ exp_{\sigma(x)})\right)_{\ast}$, hence $D((F^{-1})_{\ast
}\circ Exp_{\sigma})$ is smooth and tame.

Also, \textbf{for any} $g\in\widetilde{W}$, $D((F^{-1})_{\ast}\circ
Exp_{\sigma})(g,\cdot):\Gamma(\sigma^{\ast}(VTE))\rightarrow\Gamma(f^{\ast
}(VTE))$ is invertible and the family of inverses $V\left(  (F^{-1})_{\ast
}\circ Exp_{\sigma}\right)  :\widetilde{W}\times\Gamma(f^{\ast}
(VTE))\rightarrow\Gamma(\sigma^{\ast}(VTE))$ is given by $\left(\coprod_{x\in M}G_x\right)_{\ast}$, which is also smooth and tame.

The conclusion then follows from the Nash-Moser inverse function theorem.
\end{proof}
\end{theorem}

We want to end this Appendix by looking at one more example. This example is
actually quite trivial and can be found in \cite{ACM89} p.545 and in \cite{KR86}
p.20 (in a different form), but we want to emphasize that it can be viewed as
a very special case of the $L^2$ metric.

\begin{example}
Let $P$ be a principal $G$-bundle over an oriented closed manifold $M$. The
space of connection on $P$, denoted by $Conn(P)$, carries an affine
structure with the translation vector space $\Omega^1(M,AdP)$, in this way
$Conn(P)$ is a Fr\'echet manifold with its tangent space at each point
canonically identified with $\Omega^1(M,AdP)$. Any choice of a Riemannian
metric $g$ on $M$ and an $Ad$-invariant inner product on the Lie algebra of
$G$ gives the vector bundle $T^{\ast}M\otimes AdP$ a Riemannian structure
$\left\langle\cdot,\cdot\right\rangle$, and further gives an $L^2$ inner
product $(\cdot,\cdot)_{L^2}$ on $\Omega^1(M,AdP)$. For any fixed
connection $A\in Conn(P)$, we just define the inner product on $T_{A}
Conn(P)\cong\Omega^1(M,AdP)$ by:
\begin{equation}
(u,v)_{A}\triangleq(u,v)_{L^2}=\int_M\left\langle u(x),v(x)\right\rangle\mu_g(x).
\label{L2-con}
\end{equation}

The above formula defines a flat metric on $Conn(P)$, and this metric is an
example of $L^2$ metric in our sense (each fiber of $T^{\ast}M\otimes AdP$
is endowed with the flat Riemannian structure determined by the inner product
on it). The exponential map at $A$ is just $\Omega^1(M,AdP)\rightarrow
Conn(P)$ sending $v$ to $A+v$, which is obviously a diffeomorphism. This
trivial exponential map, together with the fact that the \eqref{L2-con} is
invariant under the gauge transformation group, has been used to construct the slice for the
gauge group acting on the space of connections, see \cite{ACM89} and
\cite{KR86}.
\end{example}
\begin{remark}
Besides the above-mentioned example, it is also well-known that for a closed
oriented manifold, the action of the diffeomorphism group on the space of
Riemannian metrics also admits a slice (\cite{E70}, \cite{Sub85}). Does the
same thing happens in our case, that is, for the gauge group $GL(E)$ acting on
the space of Hermitian metrics? (Note that we have already proved in
\eqref{acting} that our $L^2$ metric is invariant under $GL(E)$.) The author
thinks in this case the slice may not exist, for $E$ is not compact so we
cannot conclude that the action of $GL(E)$ on $\Gamma(Herm^+(E))$ is proper;
as a result, the usual routine proof falls apart.
\end{remark}

\bibliographystyle{abbrv}
\bibliography{ref}

\end{document}